%%%%%%%%%%%%%%%%%%%%%%%%%%%%%%%
%
% This is abgrp-ccr4.tex
%
% Last update by K-H, 17.11.08
%
%%%%%%%%%%%%%%%%%%%%%%%%%%%%%%%

\input amssym.def
\input amssym.tex
\openup 5pt
\pageno=1

%Macros for the file

% Achtung: Bei Verwendung von jltmac werden die Fontfamilien masm und
% msbm mehrfach aufgerufen. Das kann zu Kapazitaetsschwierigkeiten
% fuehren. KH 27.2.95

\def\item#1{\vskip1.3pt\hang\textindent {\rm #1}}% THIS REPLACES KNUTH'S DEF'N

                                  % THIS REPLACES KNUTH'S DEF'N

\newskip\litemindent
\litemindent=0.7cm  %%%%% amount of indentation for \litem
\def\Litem#1#2{\par\noindent\hangindent#1\litemindent
\hbox to #1\litemindent{\hfill\hbox to \litemindent
{\ninerm #2 \hfill}}\ignorespaces}
\def\litem{\Litem1}

\tolerance=300
\pretolerance=200
\hfuzz=1pt
\vfuzz=1pt

%\magnification=\magstep1    % am 19.6.95 auskommentiert
\hoffset=0in
%\voffset=0.2in

%\baselineskip
\hsize=5.8 true in
\vsize=9 true in
%\baselineskip
\parindent=25pt
\mathsurround=1pt
\parskip=1pt plus .25pt minus .25pt
\normallineskiplimit=.99pt

\countdef\revised=100
\mathchardef\emptyset="001F % THIS REPLACES KNUTH'S DEFINITION
\chardef\ss="19
\def\3{\ss}
\def\anf{$\lower1.2ex\hbox{"}$}
\def\frac#1#2{{#1 \over #2}}
\def\>{>\!\!>}
\def\<{<\!\!<}

\def\into{\hookrightarrow}

\def\ssssarr{\hbox to 15pt{\rightarrowfill}}
\def\sssarr{\hbox to 20pt{\rightarrowfill}}
\def\ssarr{\hbox to 30pt{\rightarrowfill}}
\def\sarr{\hbox to 40pt{\rightarrowfill}}
\def\arr{\hbox to 60pt{\rightarrowfill}}
\def\larr{\hbox to 60pt{\leftarrowfill}}
\def\Arr{\hbox to 80pt{\rightarrowfill}}

\def\cl{\mathop{\rm cl}\nolimits}

\def\Hom{\mathop{\rm Hom}\nolimits}%
\def\id{\mathop{\rm id}\nolimits} % USED FOR IDENTITY FUNCTION
\def\im{\mathop{\rm im}\nolimits}

% USED FOR IMAGINARY PART OF COMPLEX NUMBERS

%

% USED FOR REAL PART OF COMPLEX NUMBERS

%\def\reach{\mathop{\rm reach}\nolimits}
%\def\SAut{\mathop{\rm SAut}\nolimits}
%\def\SEnd{\mathop{\rm SEnd}\nolimits}

\def\span{\mathop{\rm span}\nolimits}

%\def\SS{\mathop{\rm S}\nolimits}

% USED FOR TRACE OF MATRIX

\def\trile{\trianglelefteq}

\def\0{{\bf 0}}
\def\1{{\bf 1}}

\def\b{{\frak b}}

\def\j{{\frak j}}

\def\s{{\frak s}}

\def\C{{{\Bbb C}{\mskip+1mu}}} % a \C with better spacing
 % a \D with better spacing
 % a \K with better spacing
 % a \H with better spacing

\def\R{{\Bbb R}}
\def\Z{{\Bbb Z}}
\def\N{{\Bbb N}}

\def\E{{\Bbb E}}

\def\P{{\Bbb P}}

\def\T{{\Bbb T}}

\def\:{\colon}  %8.5.92
\def\.{{\cdot}}
\def\|{\Vert}
\def\bsk{\bigskip}

\def\giantskip{\vskip2\bigskipamount}
\def\gsk{\giantskip}

\def\msk{\medskip}

\def \res {\!\mid\!\!}

\def\bbr{\bigbreak}
\def\giantbreak{\par \ifdim\lastskip<2\bigskipamount \removelastskip
         \penalty-400 \giantskip\fi}

\def\nin{\noindent}
\def\cen{\centerline}
\def\pagebreak{\vskip 0pt plus 0.0001fil\break}
\def\linebreak{\break}

\def\eps{\varepsilon}
\def\epsilon{\varepsilon}

\def\nin{\noindent}

\def\pder#1,#2,#3 { {\partial #1 \over \partial #2}(#3)}
\def\pde#1,#2 { {\partial #1 \over \partial #2}}
\def\phi{\varphi}

% Besser \ltimes und \rtimes aus dem AMS-Symbols verwenden.
%\def\sdir#1{\hbox{$\mathrel\times{\hskip -4.6pt
%            {\vrule height 4.7 pt depth .5 pt}}\hskip 2pt_{#1}$}}

\def\subeq{\subseteq}

\def\Rarrow{\Rightarrow}

\def\tilde{\widetilde}

\font\ninerm=cmr9
\font\eightrm=cmr8

\font\eightbf=cmbx8

% SANS SERIF 10 POINT
 %SANS SERIF 10 POINT ITALIC

\font\smc=cmcsc10
%\font\smc8=cmcsc8
 %SLANTED TYPEWRITER 10 POINT
 %BOLD FACE MATH SYMBOLS 10 POINT
 %DUNHILL STYLE 10 POINT
 %SAN SERIF BOLD EXTENDED 10 POINT
 %USED FOR TITLES
 %USED FOR TITLES
\font\bfone=cmbx10 scaled\magstep1 %BOLDFACE AT MAGSTEP 1
\font\bftwo=cmbx10 scaled\magstep2 %BOLDFACE AT MAGSTEP 2
 %BOLDFACE AT MAGSTEP 3

\def\qed{{\unskip\nobreak\hfil\penalty50\hskip .001pt \hbox{}\nobreak\hfil
          \vrule height 1.2ex width 1.1ex depth -.1ex
           \parfillskip=0pt\finalhyphendemerits=0\medbreak}\rm}
%This is the end-of-proof sign.
%Not to be used in display mode.
%If you want to conclude a proof
%at the end of a line is display mode use

\def\qeddis{\eqno{\vrule height 1.2ex width 1.1ex depth -.1ex} $$
                   \medbreak\rm}
%BUT OMIT $$---the macro will write that

\def\Lemma #1. {\bigbreak\vskip-\parskip\noindent{\bf Lemma #1.}\quad\it}

\def\Sublemma #1. {\bigbreak\vskip-\parskip\noindent{\bf Sublemma #1.}\quad\it}

\def\Proposition #1. {\bigbreak\vskip-\parskip\noindent{\bf Proposition #1.}
\quad\it}

\def\Corollary #1. {\bigbreak\vskip-\parskip\nin{\bf Corollary #1.}
\quad\it}

\def\Theorem #1. {\bigbreak\vskip-\parskip\noindent{\bf Theorem #1.}
\quad\it}

\def\Definition #1. {\rm\bigbreak\vskip-\parskip\noindent
{\bf Definition #1.}
\quad}

\def\Remark #1. {\rm\bigbreak\vskip-\parskip\noindent{\bf Remark #1.}\quad}

\def\Example #1. {\rm\bigbreak\vskip-\parskip\noindent{\bf Example #1.}\quad}
\def\Examples #1. {\rm\bigbreak\vskip-\parskip\noindent{\bf Examples #1.}\quad}

\def\Problems #1. {\bigbreak\vskip-\parskip\noindent{\bf Problems #1.}\quad}
\def\Problem #1. {\bigbreak\vskip-\parskip\noindent{\bf Problem #1.}\quad}
\def\Exercise #1. {\bigbreak\vskip-\parskip\noindent{\bf Exercise #1.}\quad}

\def\Conjecture #1. {\bigbreak\vskip-\parskip\noindent{\bf Conjecture #1.}\quad}

\def\Proof#1.{\rm\par\ifdim\lastskip<\bigskipamount\removelastskip\fi\smallskip
            \noindent {\bf Proof.}\quad}

\def\Axiom #1. {\bigbreak\vskip-\parskip\noindent{\bf Axiom #1.}\quad\it}

\def\Satz #1. {\bigbreak\vskip-\parskip\noindent{\bf Satz #1.}\quad\it}

\def\Korollar #1. {\bbr\vskip-\parskip\nin{\bf Korollar #1.} \quad\it}

\def\Folgerung #1. {\bbr\vskip-\parskip\nin{\bf Folgerung #1.} \quad\it}

\def\Folgerungen #1. {\bbr\vskip-\parskip\nin{\bf Folgerungen #1.} \quad\it}

\def\Bemerkung #1. {\rm\bigbreak\vskip-\parskip\noindent{\bf Bemerkung #1.}
\quad}

\def\Beispiel #1. {\rm\bigbreak\vskip-\parskip\noindent{\bf Beispiel #1.}\quad}
\def\Beispiele #1. {\rm\bigbreak\vskip-\parskip\noindent{\bf Beispiele #1.}\quad}
\def\Aufgabe #1. {\rm\bigbreak\vskip-\parskip\noindent{\bf Aufgabe #1.}\quad}
\def\Aufgaben #1. {\rm\bigbreak\vskip-\parskip\noindent{\bf Aufgabe #1.}\quad}

\def\Beweis#1. {\rm\par\ifdim\lastskip<\bigskipamount\removelastskip\fi
           \smallskip\noindent {\bf Beweis.}\quad}

\nopagenumbers

\def\date{\ifcase\month\or January\or February \or March\or April\or May
\or June\or July\or August\or September\or October\or November
\or December\fi\space\number\day, \number\year}

\def\title{Title ??}
\def\author{Author ??}

\def\thanks#1{\footnote*{\eightrm#1}}

\def\rightheadline{\hfil{\eightrm\title}\hfil\tenbf\folio}
\def\leftheadline{\tenbf\folio\hfil{\eightrm\author}\hfil}
\headline={\vbox{\line{\ifodd\pageno\rightheadline\else\leftheadline\fi}}}

\def\firstheadline{}
\def\firstfootline{\cen{\rm\folio}}

\def\seite #1 {\pageno #1
               \headline={\ifnum\pageno=#1 \firstheadline
               \else\ifodd\pageno\rightheadline\else\leftheadline\fi\fi}
               \footline={\ifnum\pageno=#1 \firstfootline\else{}\fi}}

%%%THIS IS THE MACRO LEFTSPACE.TEX %%%TO THD VIA WAFRUPP
\newdimen\dimenone
 \def\checkleftspace#1#2#3#4{%DIESER MACRO STAMMT VON APPELT
 \dimenone=\pagetotal%#1=Skip vorher,#2=Font,#3=Text,#4=Skip nachher
 \advance\dimenone by -\pageshrink   %testen ob Titel noch mit Gewalt auf Seite
                                                                          %geht
 \ifdim\dimenone>\pagegoal          %nacha tua nix-- gewoehnliche Outputroutine
   \else\dimenone=\pagetotal
        \advance\dimenone by \pagestretch
        \ifdim\dimenone<\pagegoal
          \dimenone=\pagetotal
          \advance\dimenone by#1         %addieren Skip vor Ueberschrift (=#1)
          \setbox0=\vbox{#2\parskip=0pt                %#2 ist gewaehlter Font
                     \hyphenpenalty=10000
                     \rightskip=0pt plus 5em
                     \noindent#3 \vskip#4}    %#3=Ueberschrift,#4=skip nachher
        \advance\dimenone by\ht0
        \advance\dimenone by 3\baselineskip
        \ifdim\dimenone>\pagegoal\vfill\eject\fi
          \else\eject\fi\fi}

%%% OUR HEADLINE MACROS LOOK LIKE THIS USING THIS MACRO

\def\subheadline #1{\nin\bigbreak\vskip-\lastskip
      \checkleftspace{0.9cm}{\bf}{#1}{\medskipamount}
%      \checkleftspace{0.7cm}{\bf}{#1}{\medskipamount}
          \indent\vskip0.7cm\centerline{\bf #1}\medskip}
\def\subsection{\subheadline}

\def\lsubheadline #1 #2{\bigbreak\vskip-\lastskip
      \checkleftspace{0.9cm}{\bf}{#1}{\bigskipamount}
%     \checkleftspace{0.7cm}{\bf}{#1}{\bigskipamount}
         \vbox{\vskip0.7cm}\cen{\bf #1}\msk \cen{\bf #2}\bsk}

\def\sectionheadline #1{\bigbreak\vskip-\lastskip
      \checkleftspace{1.1cm}{\bf}{#1}{\bigskipamount}
         \vbox{\vskip1.1cm}\cen{\bfone #1}\bsk}
\def\section{\sectionheadline}

\def\lsectionheadline #1 #2{\bigbreak\vskip-\lastskip
      \checkleftspace{1.1cm}{\bf}{#1}{\bigskipamount}
         \vbox{\vskip1.1cm}\cen{\bfone #1}\msk \cen{\bfone #2}\bsk}

\def\lchapterheadline #1 #2{\bigbreak\vskip-\lastskip\indent\vskip3cm
                       \cen{\bftwo #1} \msk \cen{\bftwo #2} \gsk}
\def\llsectionheadline #1 #2 #3{\bigbreak\vskip-\lastskip\indent\vskip1.8cm
\cen{\bfone #1} \msk \cen{\bfone #2} \msk \cen{\bfone #3} \nobreak\bsk\nobreak}

%\def\[#1 #2\par{\hbox{\vtop{\hsize = 2.5 true cm \nin [#1]\hfill}
%\vtop{\hsize = 12.0 true cm \nin #2\penalty10000\llap.}}
%\vbox{\vskip.3\baselineskip}}
% parameters for hsize are percentages !!! f.e. 0.2 + 0.8 = 1.0

\newtoks\literat
\def\[#1 #2\par{\literat={#2\unskip.}%
\hbox{\vtop{\hsize=.15\hsize\nin [#1]\hfill}
\vtop{\hsize=.82\hsize\nin\the\literat}}\par
\vskip.3\baselineskip}

\def\references{
\sectionheadline{\bf References}
\frenchspacing

\entries\par}

\mathchardef\emptyset="001F

\def\abstract #1{{\narrower\baselineskip=10pt{\noindent
\eightbf Abstract.\quad \eightrm #1 }
\bigskip}}

\def\firstpage{\nin
{\obeylines \parindent 0pt }
\vskip2cm
\centerline{\bfone\title}
\gsk
\centerline{\bf\author}
\vskip1.5cm \rm}

\def\Box #1 { \msk\par\nin
\centerline{
\vbox{\offinterlineskip
\hrule
\hbox{\vrule\strut\hskip1ex\hfil{\smc#1}\hfill\hskip1ex}
\hrule}\vrule}\msk }

% Antidiagonale Punkte in Matrix
\def\adots{\mathinner{\mkern1mu\raise1pt\vbox{\kern7pt\hbox{.}}
                        \mkern2mu\raise4pt\hbox{.}
                        \mkern2mu\raise7pt\hbox{.}\mkern1mu}}

\def\Rep{\mathop{\rm Rep}\nolimits}

\def\ilim{\mathop{{\rm lim}}\limits_{\longrightarrow}}
\def\slim{\mathop{{\rm s-}\!\lim}}
\def\rest{\mathord{\restriction}}
\def\un{\hbox{\mathsurround=0pt${\rm 1}\!\!{\rm 1}$\mathsurround=5pt}}
\def\cl#1.{{\cal #1}}
\def\al#1.{{\cal #1}}
\def\b#1.{{\bf #1}}
\def\wt{\widetilde}

\def\chop{\hfill\break}
\def\j{\phi}
\def\CCRX{\overline{\Delta( S,\,B)}}
\def\s #1.{_{\smash{\lower2pt\hbox{\mathsurround=0pt $\scriptstyle #1$}}\mathsurround=5pt}}

% End of macros
%%%%%%%%%%%%%%%%%%%%%%%%%%%%%%%%%%%%%%%%%%%%%%%%%%%%%%%%%%%%%%%%%%%%%%%%%%

\def\title{Full regularity for a C*-algebra of the Canonical Commutation Relations}
\def\author{Hendrik Grundling, Karl-Hermann Neeb}
\def\date{1/7/08}
% F\"ur endg\"ultige Version l\"oschen
%\def\rightheadline{\tenbf\folio\hfil{\tt abgrp10.tex}\hfil\eightrm\date}
\def\sst{\scriptstyle}

\firstpage
\abstract{The Weyl algebra,- the usual C*-algebra employed to model the canonical commutation
relations (CCRs), has a well-known defect in that it has a large number of representations which are
not regular  and these cannot model physical fields. Here, we construct explicitly a
C*-algebra which can reproduce the CCRs of a
countably dimensional symplectic space ${\scriptstyle (S,B)}$ and such that its representation set is exactly
the full set of regular representations of the CCRs.
This construction uses Blackadar's version of infinite tensor products of nonunital C*-algebras,
and it produces a ``host algebra'' (i.e. a generalised group algebra, explained below) for the
${\sst\sigma}\hbox{--representation}$ theory of the abelian group ${\sst S}$
where ${\sst\sigma(\cdot,\cdot):=e^{iB(\cdot,\cdot)/2}}.$
As an easy application, it then follows that for every regular representation
of ${\sst \CCRX}$ on a separable Hilbert space, there is a
direct integral decomposition of it into irreducible regular representations
(a known result).
\hfill\break
}

\subheadline{Introduction}

In the description of quantum systems one typically
deals  with a set of operators satisfying canonical
commutation relations. This means that there is
 a real linear map $\j$  from a given symplectic space
${(S,\,B)}$ to a linear space of selfadjoint operators
on some common dense invariant core $\al D.$ in a Hilbert space $\al H.,$
satisfying the relations
$$
\big[\j(f),\,\j(g)\big]=iB(f,\, g) \, {\un}, \quad \j(f)^* =\j(f)\quad\hbox{on}\quad
\al D. \, .
$$
If $\{q_i,p_i\,\mid\,i\in I\}\subset S$ is a symplectic basis for $S$ i.e.,
$0=B(q_i,q_j)=B(p_i,p_j)=B(p_i,q_j)-\delta_{ij}$ then
$\j(q_i)$ and $\j(p_i)$ are interpreted as quantum mechanical
position and momentum operators. If $S$ consists of Schwartz functions
on a space-time manifold, we can take $\j$ to be a bosonic quantum field.

As is known, if ${(S,\,B)}$ is non--degenerate
then the operators $\j(f)$ cannot all be bounded, so
it is natural to go from the polynomial
algebra $\al P.$ generated by ${\{\j(f)\,\big|\,f\in S\}}$ to a C*--algebra
encoding the same algebraic information.
The obvious way to do this,
is to form suitable bounded functions of the fields $\j(f).$
Following Weyl, we
consider the C*--algebra generated by the set of unitaries
$$
\big\{\exp\big(i\j(f)\big)\,\big|\,f\in S\big\}
$$
and this C*--algebra is simple.  It can be defined abstractly
as the C*--algebra generated by a set of unitaries
${\big\{\delta_f \mid{f\in X}\big\}}$ subject to the relations
$\delta_f^*=\delta_{-f}$ and $\delta_f\delta_g=e^{-iB(f,g)/2}
\delta_{f+g}$.
This is the familiar Weyl (or CCR) algebra, often denoted $\CCRX$
(cf.~[M68]). A different C*-algebra for the CCRs was defined in [BG08] based on
the resolvents of the fields.

By its definition, $\CCRX$  has a representation in which
the unitaries $\delta_f$ can be identified with
the exponentials $e^{i \j(f)},$ and hence we can obtain the
concrete algebra $\al P.$ back from these.
Such representations $\pi:\CCRX\to\al B.(\al H.)\,,$ i.e. those for which
the one--parameter groups $\lambda\to\pi(\delta\s\lambda f.)$ are strong operator
continuous for all $f\in X$ are called {\it regular}, and states are regular
if their GNS--representations are.
Since for physical situations the quantum fields are defined as the generators
of the one--parameter groups $\lambda\to\pi(\delta\s\lambda f.),$ the
representations of interest are required to be regular.
(Note that the ray--continuity of $s\to\pi(\delta\s s.)$ implies
continuity on all finite dimensional subspaces of $S.$)
Unfortunately, $\CCRX$ has a large number of nonregular representations, and so one can object that
it is not satisfactory, since analysis of physical objects can lead to nonphysical ones, e.g.
w*-limits of regular states can be nonregular.
Nonregular representations are
interpreted as situations where the field $\j(f)$ can have ``infinite field strength''.
Whilst this is useful for some nonphysical idealizations e.g. plane waves (cf.~[AMS93]),
or for quantum constraints (cf.~[GH88]), for physical situations
one wants to exclude such representations.
The resolvent algebra of [BG08] also has nonregular representations
(although far fewer than the Weyl algebra).
Our aim here is to construct a C*-algebra for the CCRs
of a countably dimensional $(S,B)$
such that its representation space
comprises of exactly the regular representations of the CCRs, in a sense to be made precise below.
This will demonstrate that the regular representation theory of the Weyl algebra is
isomorphic to the full representation theory of a C*-algebra, and hence it is subject
to the usual structure theory for the full representation theory of C*-algebras.
The existence of such an algebra has already been shown in~[Gr97], but here
we want to obtain an explicit construction of it.

In the case that $S$ is finite dimensional, there is an immediate solution.
Regard $S$ with its addition as an Abelian group, then
$\sigma(\cdot,\cdot):=\exp[iB(\cdot,\cdot)/2]$ is a 2-cocycle of $S,$ and
$\CCRX$ is just the $\sigma\hbox{--twisted}$ group algebra of $S$ with the discrete
topology (cf.~[PR89]).
Define the C*--algebra $\cl L.$ as the C*--envelope
 of the twisted convolution
algebra of $S,$ where the latter
 consists of $L^1(S)$ equipped with the multiplication and involution:
$$f*g({ x})=\int_{S}f({ y})\, g({x}-{ y})\,\sigma({ y},\,{ x})\,d\mu({ y})\,,
\qquad\; f^*({ x})=\overline{f(-{ x})}$$
where $\mu$ is a Haar measure on $S,$ i.e. $\cl L.$ is
the $\sigma\hbox{--twisted}$ group algebra of $S\,.$
This algebra $\cl L.$ is known to be isomorphic to the compacts $\cl K.(L^2(S))$
(cf.~[Se67] and [Bla06] p206).
Then we have an embedding of $\CCRX$ into the multiplier algebra of $\cl L.,$
$\CCRX\subset M(\cl L.),$ by the action
$\delta_{ x}\cdot  f({ y})= \sigma({ x},\,{ y})\, f({ y}-{ x}).$
 The unique extensions of representations on $\cl L.$ to $\CCRX$ produces
a bijection from the representations of $\cl L.$ onto the
regular representations of $\CCRX$ and the bijection respects direct sums and takes
irreducibles to irreducibles. So $\cl L.$ is the desired C*-algebra with full
regularity.

For the case that $S$ is infinite dimensional,
since regular representations $\pi$ are characterized by requiring
the maps $s\to\pi(\delta\s s.)$ to be  continuous on all finite dimensional
subspaces of $S,$ this means that we require these maps to be strong operator
continuous w.r.t. the inductive limit topology, where the inductive limit is
the one consisting of all finite dimensional subspaces of $S$ under inclusion.
This inductive limit topology on $S$ is only a group topology w.r.t. addition in the case that
$S$ is a countably dimensional space; cf.~[Gl03]. Hence in this case the
regular representation theory of $\CCRX$ is the $\sigma\hbox{--representation}$ theory of the topological
group $S,$ but not otherwise.
% (this may explain why so many important results
% hold for the regular representations of this case, but not otherwise).
Henceforth we will always take $(S,B)$ to be countably dimensional,
equipped with
the (locally convex) inductive limit topology.
The problem now becomes
the one of how to
define a $\sigma\hbox{--twisted}$ group algebra for $S.$ The usual theory fails
in this case, since $S$ is not locally compact, hence there is no Haar measure.

We see that there is a need to generalize
the notion of a (twisted) group algebra to topological groups which are not locally compact.
Such a generalization,
called a {\it full host algebra,} has been proposed in~[Gr05]. Briefly,
it is a $C^*$-algebra $\cl A.$ which has in its multiplier algebra
$M(\cl A.)$ a homomorphism $\eta \: G \to U(M({\cal A}))$,
such that the (unique) extension of
the representation theory of $\cl A.$ to $M(\cl A.)$ pulls back
via $\eta$ to the continuous (unitary) representation theory
of $G$. There is also an analogous concept for
unitary $\sigma\hbox{--representations}$, where $\sigma$ is a
continuous $\T$-valued $2$-cocycle on $G$.
Thus, given a full host algebra $\cl A.,$
the continuous representation theory of $G$ can be analyzed on
$\cl A.$ with a large arsenal of $C^*$-algebraic tools.

Our main result in this paper is
an explicit construction of a full host algebra for the
$\sigma\hbox{--representations}$ of an infinite dimensional
topological linear space $S,$ regarded as a group where
$S$ will be a countably dimensional symplectic space
with symplectic form $B,$ equipped with
the (locally convex) inductive limit topology.
We demonstrate the usefulness of this construction by proving that
for every regular representation
of $\CCRX$ on a separable Hilbert space, there is a
direct integral decomposition of it into irreducible regular representations.
This last result is already known by different means (cf.~[He71,~Sch90]).

This paper is structured as follows. In Section~I we state the notation and
definitions  necessary for the subsequent material,
and in Section~II we discuss existence and uniqueness issues for host algebras.
In Section~III we construct the
host algebra for the pair $(S,\sigma)$ mentioned above,
do the direct integral decomposition mentioned,
and in the appendix
we add general results concerning host algebras and the strict topology
which are required for our proofs. These results are of independent
interest for the general structure theory of host algebras.
The reader in a hurry can skip Section~II.

\subheadline{I. Definitions and notation}

%{\bf Preliminaries and notation:}
We will need the following notation and concepts for our main results.
\item{$\bullet$}
In the following, we write $M({\cal A})$ for the
multiplier algebra of a $C^*$-algebra ${\cal A}$ and, if ${\cal A}$ has a unit, $U({\cal A})$ for its
unitary group.
We have an injective morphism of $C^*$-algebras
$\iota_{\cal A} \: {\cal A} \to M({\cal A})$ and will just denote ${\cal A}$ for its
image in $M({\cal A})$. Then ${\cal A}$ is dense
in  $M({\cal A})$ with respect to
the {\it strict topology},
which is the locally convex topology defined by the seminorms
$$ p_a(m) := \|m \cdot a\| + \|a \cdot m\|,
\qquad a\in {\cal A},\; m\in M({\cal A})$$
(cf.\ [Wo95]).
\item{$\bullet$}
For a complex Hilbert space ${\cal H}$, we write $\Rep({\cal A},{\cal H})$ for the
set of non-degenerate representations of ${\cal A}$ on ${\cal H}$.
Note that the collection $\Rep {\cal A}$ of all non-degenerate
representations of ${\cal A}$ is not a set, but a (proper) class
in the sense of von Neumann--Bernays--G\"odel set theory, cf.~[TZ75], and in
this framework we can consistently manipulate the object $\Rep {\cal A}.$
However, to avoid set--theoretical subtleties, we will express our results below
concretely, i.e., in terms of $\Rep({\cal A},{\cal H})$ for given Hilbert spaces
$\cl H..$
We have an injection
$$ \Rep({\cal A}, {\cal H}) \into \Rep(M({\cal A}),{\cal H}), \quad \pi \mapsto \tilde\pi
\quad \hbox{ with } \quad \tilde\pi \circ \iota_{\cal A} = \pi, $$
which identifies the non-degenerate representation $\pi$ of
${\cal A}$ with that representation $\tilde\pi$ of its multiplier algebra which
extends $\pi$ and is
continuous with respect to the  strict topology on $M({\cal A})$ and
the topology of pointwise convergence on $B({\cal H})$.
\item{$\bullet$}
For  topological groups $G$ and $H$ we write
$\Hom(G,H)$ for the set of continuous group homomorphisms $G \to H$.
We also write
$\Rep(G,{\cal H})$ for the set of all (strong operator)
continuous unitary representations of
$G$ on ${\cal H}\,.$ %and $\Rep G$ for the (proper) class
% of continuous unitary representations of $G$.
Endowing $U({\cal H})$ with the strong operator topology
turns it into a topological group, denoted $U({\cal H})_s$,
so that $\Rep(G,{\cal H}) = \Hom(G,U({\cal H})_s)$.
\item{$\bullet$}
 Let $\T \subeq \C^\times$ denote the unit circle, viewed as a multiplicative subgroup
and $\sigma \: G \times  G\to \T$ be a continuous $2$-cocycle, i.e.,
$$ \sigma(\1,x) = \sigma(x,\1) = 1, \quad
\sigma(x,y) \sigma(xy,z)  = \sigma(x,yz)\sigma(y,z)
\quad \hbox{ for } \quad x,y,z \in G. $$
We then form the topological group
$$ G_\sigma := \T \times G, \quad
(t,g)(t',g') := (tt' \sigma(g,g'), gg') $$
and note that the projection $q \: G_\sigma \to G$ defines a central
extension of $G$ by $\T$.
A continuous unitary representation $(\pi, {\cal H})$
of $G_\sigma$ is called a {\it $\sigma$-representation
of $G$} if $\pi(t,\1) = t \1$ holds for each $t \in \T$. Then
$$ G \to U({\cal H}), \quad g \mapsto \pi(1,g) $$
is continuous with respect to the strong operator topology, but
$$ \pi(1,g)\pi(1,g') = \sigma(g,g') \pi(1,gg') \quad \hbox{ for } \quad g,g' \in G. $$
We write $\Rep((G,\sigma), {\cal H})$ for the set of all continuous
$\sigma$-representations of $G$ on ${\cal H}$.

%
% \subheadline{I. Definitions and main results}
%

\Definition I.1. Let  $G$ be a topological group
and $\sigma\: G \times  G\to \T$ a continuous 2-cocycle.\chop
A {\it host algebra for the pair $(G,\sigma)$} is a pair
$({\cal L}, \eta)$ where ${\cal L}$ is a $C^*$-algebra
and $\eta \: G_\sigma \to U(M({\cal L}))$ is a homomorphism
such that for each complex Hilbert space
${\cal H}$ the corresponding map
$$ \eta^* \: \Rep({\cal L},{\cal H}) \to \Rep((G,\sigma), {\cal H}), \quad
\pi \mapsto \tilde\pi \circ \eta $$
is injective. We then write $\Rep(G,{\cal H})_\eta \subeq \Rep(G,{\cal H})$ for the range
of $\eta^*$.
We say that $(G,\sigma)$ has a {\it full host algebra} if it has a host algebra
for which $\eta^*$ is surjective for each Hilbert space~${\cal H}$.

In the case that $\sigma=1$, we simply  speak of a {\it host algebra for $G$}.
In this case, $G_\sigma = G \times \T$ is a direct product, so that
a host algebra for $G$ is a pair
$({\cal L}, \eta)$, where $\eta \: G \to U(M({\cal L}))$ is a homomorphism
into the unitary group of $M({\cal L})$ such that for each complex Hilbert space
${\cal H}$ the corresponding map
$$ \eta^* \: \Rep({\cal L},{\cal H}) \to \Rep(G,{\cal H}), \quad
\pi \mapsto \tilde\pi \circ \eta $$
is injective. We then write $\Rep(G,{\cal H})_\eta \subeq \Rep(G,{\cal H})$ for the range
of $\eta^*$.
We say that $G$ has a {\it full host algebra} if it has a host algebra
for which $\eta^*$ is surjective for each Hilbert space~${\cal H}$.

Note that by the universal property of (twisted) group algebras, the homomorphism
$\eta \: G_\sigma \to U(M({\cal L}))$  extends uniquely to the $\sigma\hbox{--twisted}$
group algebra of $G$ with the discrete topology, i.e., we have a *-homomorphism
$\eta \: C^*_\sigma(G_d)\to U(M({\cal L}))$ (still denoted by $\eta$).
\qed

\item{\bf Remark:} (1)
It is well known that for each locally compact group $G$,
the group $C^*$-algebra $C^*(G)$, and
the natural map $\eta_G \: G \to M(C^*(G))$ provide a
full host algebra ([Dix64, Sect.~13.9]) and for each pair $(G,\sigma)$, where
$G$ is locally compact, the corresponding twisted group
$C^*$-algebra $C^*(G,\sigma)$, which is isomorphic to an ideal of $C^*(G_\sigma)$,
is a full host algebra for the pair $(G,\sigma)$. This is most easily seen
by decomposition of representations of $
G_\sigma$ into isotypic summands with respect to the
action of the central subgroup $\T \times \{\1\}$ (apply [BS70], [PR89]
with ${\cal L}=\C$). The map $\eta_G: G\to M(C^*(G,\sigma))$
is continuous w.r.t. the strict topology of
$M(C^*(G,\sigma))$.\footnote{$^1$}{This is an easy consequence of the fact that
$\im(\eta_G)$ is bounded and that the action on the corresponding $L^1$-algebra
is continuous.}
\item{(2)}
Note that the map $\eta^*$ preserves direct sums, unitary conjugation, subrepresentations,
and for full host algebras, irreducibility (cf.~[Gr05]) so that this notion of isomorphism between
$\Rep((G,\sigma), {\cal H})$ and $\Rep({\cal L},{\cal H})$ involves strong structural correspondences.
\item{(3)}
Whilst the concept of a host algebra is a natural extension of the concept of a group
C*-algebra (and it easily generalizes to other algebraic objects cf.~[Gr05]), it has so far had
a troubled history. It was first used in [Gr97], though not under this name.
There, the existence of host algebras was proven for groups which are inductive
limits of locally compact groups, though the proof was not constructive enough
to allow much further structural analysis of these host algebras.
Then in [Gr05] the concept was generalised to algebraic objects other than topological groups,
and a general existence and uniqueness theorem was given, though unfortunately
this turned out to be wrong (see the erratum, and the counterexample below).
Since then, host algebras have been constructed in [Ne08] for complex semigroups.
Our aim in Section III is to provide an explicit, and more useful construction
of a host algebra (than [Gr97])
for the regular representations of the canonical commutation relations.

\subheadline{II. Existence and Uniqueness issues.}

For general topological groups, there are serious existence and uniqueness
questions for their host algebras (as mentioned, the
existence and uniqueness theorem in~[Gr05] is  wrong).
From the structural ``isomorphism" between the $\sigma\hbox{-representation}$ theory of $G$
and the representation theory of its full host noted above,
it becomes easy to find examples of
topological groups without full host algebras.
For instance in Example~5.2 of~[Pe97] is an abelian topological group
with a faithful continuous unitary
representation, but no continuous irreducible
representations. Hence this group cannot have a host algebra, whether
full or not. In [GN01] it is shown in particular for any
non-atomic measure space $(X,\mu)$, such as the unit interval
$[0,1]$ with Lebesgue meaure, the unitary group
of the $W^*$-algebra $L^\infty(X,\mu)$, endowed with the weak topology,
has no non-trivial continuous characters, hence no non-zero host algebra.

It is therefore an important open problem to characterize those pairs
 $(G,\sigma)$ for which full host algebras exist.

Concerning the issue of uniqueness, the following simple counterexample
shows that if a host algebra exists, then it need not be unique.
Let $G := \Z.$ Then its
character group is $\widehat{G}\cong \T,$ which is a compact group with respect to the topology
of pointwise convergence. Since $G$ is locally compact, $C^*(G) \cong C(\T)$
is a full host algebra for $G.$ Let ${\cal L} := C_0([0, 1))$ and define a homomorphism
$\eta : \Z \to U(M(C_0([0, 1)))) \cong C([0, 1),\T)$ by $\eta(n)(x) := e^{2\pi inx}.$ Then
$\eta(1) : [0, 1)\to \T$ is a continuous bijection, which implies that
$\eta(\Z)$ separates the points, hence  by Lemma~A.1 below
the $C^*$-algebra generated by this set
is strictly dense in $M({\cal L})$.
Since the unique extensions of representations of ${\cal L}$ to $M({\cal L})$
are continuous in the strict topology, it follows that $\eta^*$ is injective.
Further, $\Z$ is discrete, so that continuity of representations $\eta^*\pi$
is trivially satisfied, and thus ${({\cal L},\,\eta)}$ is a host algebra. This host algebra
is full because the representations of $\Z$ are in one-to-one correspondence
with Borel spectral measures on $\T$ and $\eta(1)$ is a Borel isomorphism. Note in
particular that this full host algebra ${\cal L}\not\cong C^*(G) $ is not unital,
although $G$ is a discrete
group.

This issue also needs further analysis, e.g. one needs to find what structural
properties are shared by host algebras for the same pair  $(G,\sigma)$,
and to explore the properties of the set of host algebras.
In the appendix we list more host algebra properties, e.g. those relating to products
and homomorphisms of groups.

\subheadline{III. A construction of a full host algebra for $(S,\,\sigma)\,.$}

\noindent
Here we want to present an example of a host algebra for an infinite-dimensional
group. Let ${(S,B)}$ be a countably dimensional (nondegenerate) symplectic space.
Then by Lemma~A.8 we know that there is a complex
structure and a hermitian inner product
$(\cdot,\cdot)$ on $S$ such that $B(v,\,w)={\rm Im}\,(v,\,w)$ for all $v,\,w\in S$.
Moreover, w.r.t. the inner product $(\cdot,\cdot),$ $S$ has an orthonormal basis
 $(e_n)_{n \in \N}$.
We consider $S \cong \C^{(\N)}$
as an inductive limit of the subspaces
$S_n := \span \{ e_1,\ldots, e_n\}$ and
endow it with the inductive limit topology, which turns it into an abelian topological
group with respect to addition (which is only true for countably dimensional spaces; cf.\ [Gl03]).
Moreover, the symplectic form
$B(v,\,w)={\rm Im}\,(v,\,w)$ defines a group two-cocycle
$\sigma(v,\,w):=\exp[iB(v,w)/2]$
on $S$.
Let $S_\sigma$ denote the corresponding central extension of $S$ by $\T$
(cf. above Definition~I.1). In the rest of this section we will prove that:

\Theorem III.1. The pair $(S,\sigma)$ has a full host algebra.
\qed

% Denote by $\cal R$ the class of  % those (continuous) representations $\pi$ of
% $S_\sigma$ such that $\pi(0,t)=t\un$, i.e., they map the central group
% $(0,\T)$ of $S_\sigma$ identically to the scalars $\T\un$.
% These representations of $S_\sigma$ are of course just
%  $\sigma\hbox{--representations}$ of $S$.

Recall that $\cl A.:=\CCRX$ is the discrete twisted $\sigma\hbox{--group}$ algebra
of $S$, i.e., it is the unique (simple) $C^*$-algebra generated by a collection
of unitaries ${\big\{\delta_s\;\big|\;s\in S\big\}}$ satisfying the
(Weyl) relations $\delta_{s_1}\delta_{s_2}=\sigma(s_1,s_2)\,\delta_{s_1+s_2}$
([BR97, Th.~5.2.8]). Let
$$\cl R.(\cl H.):=\big\{\pi\in\Rep(\cl A.,{\cal H})\,\mid\,
t\in\R\to\pi(\delta\s tx.)\;\hbox{is strong operator continuous}\;\forall x\in S\big\}
$$
denote the set of regular representations on the Hilbert space $\cl H..$
Through the identification $\pi(s):=\pi(\delta_{s}),$ $\cl R.(\cl H.)$ corresponds
exactly with the $\sigma$-representations of  $S$ on $\cl H.,$ i.e.,
with ${\Rep((S,\sigma), \cl H.)\,.}$

\Lemma III.2.  With the notation above, we have ${\cal A}=
\mathop{\bigotimes}\limits_{n=1}^\infty{\cal A}_n$ with the spatial (minimal) tensor norms,
where ${\cal A}_n:= {{\rm C}^*\big\{\delta_{ze_n}\,\big|\,z\in \C\big\}}$.

\Proof.
This follows directly from Proposition 11.4.3 of Kadison and Ringrose~[KR83],
we only need to verify that its conditions hold in the present context.
For this, observe that  ${\cal A}= {{\rm C}^*\big\{
\mathop{\bigcup}\limits_{n=1}^\infty{\cal A}_n\big\}}$, $\un\in{\cal A}_n$,
${\big[{\cal A}_n,\,{\cal A}_m\big]}=\{0\}$ when $n\not= m$. Moreover, the linear maps
$$\psi_k:{\cal A}_1\otimes\cdots\otimes{\cal A}_k\to {\cal A}\qquad
\hbox{defined by}\qquad
\psi_k\big(A_1\otimes\cdots\otimes A_k):= A_1A_2\cdots A_k$$
are *--monomorphisms because each image subalgebra ${{\rm C}^*\big\{
\mathop{\bigcup}\limits_{n=1}^k{\cal A}_n\big\}}$ is the unique $C^*$-algebra generated
by the unitaries ${\big\{\delta_{ze_i}\;\big|\;z\in\C,\;i=1,\ldots,k\big\}}$,
and this is also true for ${\cal A}_1\otimes\cdots\otimes{\cal A}_k$.
This is enough to apply the proposition~loc.~cit.
\qed

Observe that each  ${\cal A}_n$ is just the discrete
$\sigma$-group algebra of the subgroup $\C e_n\subset S$, and as the latter is locally
compact, we can construct its $\sigma\hbox{--twisted}$ group algebra which
we denote by ${\cal L}_n$ (recall that ${\cal L}_n$ is just the enveloping
$C^*$-algebra of $L^1(\C),$ equipped with $\sigma\hbox{-twisted}$ convolution).
It is well-known that  ${\cal L}_n\cong{\cal K}(L^2(\R))$ (cf. Segal~[Se67]).
Note that for each finite subset $F\subset\N,$ the algebra
$\mathop{\bigotimes}\limits_{n\in F}{\cal L}_n\cong{\cal K}(
\mathop{\bigotimes}\limits_{n\in F}L^2(\R))\cong
{\cal K}(L^2(\R^F))$ is a host algebra for the regular
representations of $\mathop{\bigotimes}\limits_{n\in F}{\cal A}_n
={{\rm C}^*\big\{\delta_{ze_n}\,\big|\,z\in \C,\;n\in F\big\}}$,
i.e., for the $\sigma\hbox{--representations}$ of
$\span\{e_n\,|\,n\in F\}\subset S$.

It is natural to try some infinite tensor product
$\mathop{\bigotimes}\limits_{n=1}^\infty{\cal L}_n$
for a host algebra, but because the algebras ${\cal L}_n$
are non-unital, the definition of the infinite tensor product
needs some care~[Bla77]. For each $n\in\N$, choose a nonzero projection
$P_n\in{\cal L}_n\cong{\cal K}({\cal H})$ and define $C^*$-embeddings
$$\Psi_{\ell k}:{\cal L}^{(k)}\to
{\cal L}^{(\ell)}\qquad\hbox{by}\qquad
\Psi_{\ell k}(A_1\otimes\cdots\otimes A_k):=
A_1\otimes\cdots\otimes A_k\otimes P_{k+1}\otimes\cdots\otimes P_\ell, $$
where $k<\ell$ and ${\cal L}^{(k)}:=\mathop{\bigotimes}\limits_{n=1}^k{\cal L}_n$.
 Then the inductive limit makes sense, so we define
$$ {\cal L} := \bigotimes_{n=1}^\infty{\cal L}_n:=\ilim\big\{
{\cal L}^{(n)},\,\Psi_{\ell k}\big\}$$
and write $\Psi_k \: {\cal L}^{(k)} \to {\cal L}$ for the corresponding embeddings,
satisfying $\Psi_k \circ \Psi_{kj} = \Psi_j$ for $j \leq k$.
Since each ${\cal L}_n$ is simple, so are the finite tensor products
${\cal L}^{(k)}$ ([WO93], Prop.~T.6.25),
and as inductive limits of simple $C^*$-algebras
are simple ([KR83], Prop.~11.4.2), so is $\cl L.$.
It is also clear that $\cl L.$ is
separable.

Since
$\Psi_{k+n,k}(L_k)=L_k\otimes P_{k+1}\otimes\cdots\otimes P_{k+n},$
where $L_k\in{\cal L}^{(k)}$, this means that we can consider ${\cal L}$
to be built up out of elementary tensors of the form
$$\Psi_k(L_1 \otimes \cdots \otimes L_k)=
L_1\otimes L_2\otimes\cdots\otimes L_k\otimes P_{k+1}\otimes P_{k+2}\otimes\cdots\;,
\qquad\hbox{where}\qquad L_i\in{\cal L}_i\;, \eqno(4.1) $$
i.e., eventually they are of the form $\cdots\otimes P_{k}\otimes P_{k+1}\otimes\cdots$.
We will use this picture below, and generally will not indicate the maps $\Psi_k\,.$

\Lemma III.3.
\item{(i)} With respect to componentwise multiplication, we have an inclusion
$${\cal A}=
\bigotimes_{n=1}^\infty{\cal A}_n \subset M({\cal L})=
M\big(\bigotimes_{n=1}^\infty{\cal L}_n\big)\;.$$
\item{(ii)} There is a natural embedding $\iota_n \: M({\cal L}^{(n)})\into M({\cal L})$.
This is a topological embedding on each bounded subset of
$M({\cal L}^{(n)})$. Moreover, ${\cal L}^{(n)}$
% (resp. ${\cal A}^{(n)}$)
is dense in $M({\cal L}^{(n)})$ with respect to the
restriction of the strict topology of $M({\cal L})$.

%\nin{\bf Hendrik: It would be nice to have
%that ${\cal A}^{(n)}$ is also dense in this topology. True?}
%\chop{\bf Karl-Hermann: It would be nice, but I don't know. Not obvious.}
\item{(iii)}
Let $\pi\in{\rm Rep}({\cal L},\cl H.)$, and let $\pi_n$ denote the unique
representation which it induces on ${\cal L}^{(n)}\subset M({\cal L}^{(n)})
\subset M({\cal L})$ by strict extension. Then
$$\pi(L_1\otimes L_2\otimes\cdots)=\slim_{n\to\infty}\pi_n(L_1\otimes\cdots\otimes L_n)$$
for all $L_1\otimes L_2\otimes\cdots\in{\cal L}$ as in {\rm(4.1)}.

\Proof. (i) For each $k$ we obtain a homomorphism $\Theta_k:
\mathop{\bigotimes}\limits_{n=1}^k{\cal A}_n\to M({\cal L})$
by componentwise multiplication in the first $k$ entries
of $\cal L$, leaving all
entries further up invariant. By simplicity of its domain,
each $\Theta_k$ is a monomorphism.
{}From $\Theta_k\big(\mathop{\bigotimes}\limits_{n=1}^k
{\cal A}_n\big)\subset M({\cal L})$ for each $k \in \N$,
we obtain all the generating unitaries $\delta_s$ in $M({\cal L})$,
then they generate $\cal A$ in $M({\cal L})$ by uniqueness of
the $C^*$-algebra of the canonical commutation relations.

(ii) Now ${\cal L}={\cal L}^{(n)}\otimes{\cal B}$ for a $C^*$-algebra $\cal B$
(cf.\ Blackadar [Bla77, p.~315]), and $M({\cal L}^{(n)})$ embeds in
$M({\cal L})$ as $M({\cal L}^{(n)})\otimes\un$.
Therefore (ii) follows from Lemma~A.2.

(iii) Note that $U_n:= \Psi_n(\un)=\un\otimes\cdots\otimes\un\otimes
P_{n+1}\otimes P_{n+2}\otimes\cdots \in  M({\cal L})$
converges strictly to $\un$.
Recall that $L=L_1\otimes L_2\otimes\cdots\in {\cal L}$ as in (4.1)
is of the form
$$ A_1\otimes A_2\otimes\cdots\otimes A_k\otimes P_{k+1}\otimes P_{k+2}\otimes\cdots\;,$$
where $A_i\in{\cal L}_i$, so for $n\geq k$ we get for all $\psi\in{\cal H}_{\pi}$ that
for the strictly continuous extension $\widetilde\pi$ of $\pi$ to $M({\cal L})$:
$$\eqalignno{\big\|\widetilde\pi(L-L_1\otimes\cdots\otimes
L_n\otimes\un\otimes\un\otimes\cdots)\psi\big\|
&=\Big\|\widetilde\pi\big(L_1\otimes\cdots\otimes L_n\otimes(P_{n+1}\otimes P_{n+2}\otimes\cdots
-\un)\big)\psi\Big\|  \cr
&=\Big\|\widetilde\pi\big(L_1\otimes\cdots\otimes L_n\otimes\un\otimes\cdots\big)\cdot
\widetilde\pi(U_n-\un)\psi\Big\|\cr
&\leq C\cdot\big\|\pi(U_n-\un)\psi\big\|\to 0\cr}$$
as $n\to\infty$, where $C>0$ is chosen
such that $\|L_1\otimes\cdots\otimes L_n\|\leq C$ for all $n$, and this is possible
because $\|P_{k+1}\otimes P_{k+2}\otimes\cdots\|=1$. But this is exactly the
claim we needed to prove.
\qed

Let $\pi\in{\rm Rep}({\cal A},\cl H.)$ be regular. Observe
that $\pi$ is regular on all ${\cal A}_n$, hence there are unique
$\widehat\pi_n\in{\rm Rep}({\cal L}_n,\cl H.)$ which extend (on $\cl H.$) to
$\pi\rest{\cal A}_n$ by the host algebra property of ${\cal L}_n$.
For the distinguished projections $P_n\in{\cal L}_n$, we simplify the notation
to $\pi(P_n) := \widehat\pi_n(P_n)$. Observe that
the projections $\pi(P_j)$ all commute, and so the strong limit
$$\P_{k} := \slim_{n\to\infty}\pi(P_{k})\cdots\pi(P_n)$$
exists, and it is the projection onto the intersection of the ranges of
of all  $\pi(P_j)$, $j\geq k$. Since $\P_{k}=\pi(P_k)\,\P_{k+1}$
we have $\P_{k+1}\geq\P_{k}$ and so also $\slim\limits_{k\to\infty}\P_k
\leq\un$ exists.

We will use the notation ${\cal A}^{(n)}:=
\mathop{\bigotimes}\limits_{j=1}^n{\cal A}_j$ below.

\Proposition III.4.
Define a monomorphism $\eta:S_\sigma\to U(M({\cal L}))$ by
$\eta((s,t)):= t\delta_s\in {\cal A}\subset M({\cal L})$ (by Lemma~III.3(i)).
Then $\eta$ is continuous with respect to the strict topology on
$M({\cal L})$ and $\cal L$ is a host algebra of $(S,\sigma)$, i.e., the maps
$\eta^*:{{\rm Rep}({\cal L},\cl H.)}\to {{\rm Rep}((S,\sigma),\cl H.)}$
are injective. The range of $\eta^*$ consists of those
$\pi\in {\Rep\big((S,\sigma),\cl H.\big)}$ for which
$\slim\limits_{k\to\infty}\P_k=\un$.

\Proof. Let $\pi$ be a representation of
${\cal L}$ and $\tilde\pi$ its strictly continuous extension to $M({\cal L})$. To
see that the representation $\eta^*\tilde\pi$ of $S_\sigma$ is continuous,
we show that $\eta$ is continuous with respect to the strict topology on $M({\cal L})$.
Since $S_\sigma$ is a topological direct limit of the subgroups
$S_{m,\sigma}$, where $S_m = \span_\C \{e_1,\ldots, e_m\}$,
it suffices to show that $\eta$ is continuous on each subgroup $S_{m,\sigma}$.
Recall that the twisted group algebra
$C^*(S_m, \sigma) \cong {\cal L}^{(m)}$ is a full host algebra for $(S_m,\sigma)$
and that the corresponding strictly continuous
homomorphism $\eta_m \: S_{m,\sigma} \to M({\cal L}^{(m)})$ is
compatible with the
embedding $\iota_m \: M({\cal L}^{(m)}) \into M({\cal L})$ in the sense that
$\eta\res_{S_{m,\sigma}} = \iota_m \circ \eta_m$. Since $\iota_m$ restricts
to an embedding on the unitary group (Lemma~III.3(ii)),
the continuity of $\eta_m$ implies the continuity of $\eta$ on $S_{m,\sigma}$,
which in turn implies the continuity of $\eta$.
As a consequence, $\tilde \pi \circ \eta$ is a continuous unitary representation
of $S_{\sigma}$ for each strictly continuous representation $\tilde\pi$ of $M({\cal L})$.

To see that $\eta^*$ is injective, we have to show that two representations
$\pi_1, \pi_2$ of ${\cal L}$ for which $\eta^*\pi_1 = \eta^*\pi_2$ are equal.
If $\eta^*\pi_1 = \eta^*\pi_2$, then we obtain for each $m \in \N$ the relation
$\eta_m^*\pi_1 = \eta_m^*\pi_2$ on $S_{m,\sigma}$. This means that the corresponding
unitary representations of the group $S_{m,\sigma}$ coincide.
In view of Lemma~III.3(iii), it suffices to argue that the two
non-degenerate representations $\pi_{1,m}$ and $\pi_{2,m}$ of ${\cal L}^{(m)}$ coincide
(cf.\ Lemma~A.3 for the non-degeneracy), which in turn follows from the
host algebra property of ${\cal L}^{(m)}$ for $S_{m,\sigma}$.

To characterize the range of $\eta^*$, let $\pi\in{\rm Rep}({\cal A},\cl H.)$ be the
strictly continuous extension of a $\pi_0\in{\rm Rep}({\cal L},\cl H.)$.
Then, by Lemma~III.3(iii), it must satisfy
$$\pi_0(L_1\otimes L_2\otimes\cdots)=\slim_{n\to\infty}\pi_n(L_1\otimes\cdots\otimes L_n)$$
for all $L_1\otimes L_2\otimes\cdots\in{\cal L}$.
% Our strategy will be to show the right hand of this equation exists, and
% that it defines
% a representation $\pi_0$ on ${\cal L}$ with the right properties.
% Recall from Takesaki [Tak79, Lemma~4.1, p.~203] that
Now we have
$$\pi_n(L_1\otimes\cdots\otimes L_{n-1}\otimes P_n)
=\widetilde\pi_n(L_1\otimes\cdots
\otimes L_{n-1}\otimes\un)\,\widetilde\pi_n(\un\otimes\cdots\un\otimes P_n)$$
where $\widetilde\pi_n$ denotes the strictly continuous
extension to $M({\cal L}^{(n)})$,
 and it is obvious that
these two operators commute. From the algebra relations
${\cal A}^{(n)}\supset{\cal A}^{(n-1)}\subset M({\cal L}^{(n-1)})\subset M({\cal L}^{(n)})$,
and the host algebra properties we get that
$\widetilde\pi_n(L_1\otimes\cdots
\otimes L_{n-1}\otimes\un)=\pi_{n-1}(L_1\otimes\cdots\otimes L_{n-1})$
and $\widetilde\pi_n(\un\otimes\cdots\un\otimes P_n)=\pi(P_n)$,  so
$$\pi_n(L_1\otimes\cdots\otimes L_{n-1}\otimes P_n)
=\pi_{n-1}(L_1\otimes\cdots\otimes L_{n-1})\,\pi(P_n)\;.$$
% where we dropped the superscript `$(2)$' on the second factor.
Thus, for
$$\eqalignno{L =
L_1\otimes L_2\otimes\cdots &=
A_1\otimes A_2\otimes\cdots\otimes A_k\otimes P_{k+1}\otimes P_{k+2}\otimes\cdots
\in{\cal L},\quad\hbox{we get for $n>k$ :}  \cr
\pi_n( L_1\otimes\cdots\otimes L_n)  &= \pi_k(A_1\otimes\cdots\otimes A_k)\,
\pi(P_{k+1})\cdots\pi(P_n)\;.\cr}$$
Using the fact that the projections $\pi(P_j)$ all commute,% and so the strong limit
% $$\slim_{n\to\infty}\pi_{k+1}(P_{k+1})\cdots\pi_n(P_n)=:\P_{k+1}$$
% exists, and it is the projection onto the intersection of the ranges of
% all $\pi_n(P_n)\,,\;n>k$.
% Thus we can now define
$$\pi_0(L_1\otimes L_2\otimes\cdots)=\slim_{n\to\infty}\pi_n(L_1\otimes\cdots\otimes L_n)
=\pi_k(A_1\otimes\cdots\otimes A_k)\,\P_{k+1}\,.$$
Since $\pi_0$ is non-degenerate, and all $\pi_k\rest{\cal L}^{(k)}$ are non-degenerate,
it follows that $\slim\limits_{k\to\infty}\P_k=\un$.

Conversely, if we start from a regular representation $\pi$ of $\cal A$ which
satisfies $\slim\limits_{k\to\infty}\P_k=\un$, we will {\it define}
a representation $\pi_0$
on $\cal L$ by
$$\pi_0(L):=\pi_k(A_1\otimes\cdots\otimes A_k)\,\P_{k+1}\quad
\hbox{for}\quad L=A_1\otimes A_2\otimes\cdots\otimes A_k\otimes P_{k+1}\otimes P_{k+2}\otimes\cdots$$
where $\pi_k\in{\rm Rep}\,{\cal L}^{(k)}$ is obtained from $\pi\rest{\cal A}^{(k)}$,
using the host algebra property of ${\cal L}^{(k)}$.
To see that this can be done, note that for
$A \in {\cal L}^{(k)}$ we have
$$ \pi_k(A)\P_{k+1} = \pi_{k+1}(\Psi_{k+1,k}(A))\P_{k+2}. $$
Therefore the universal property of the direct limit algebra ${\cal L}$ implies
the existence of a % non-degenerate
representation $\pi_0$ of ${\cal L}$, satisfying
$$ \pi_0(\Psi_k(A)) = \pi_k(A)\P_{k+1} \quad \hbox{ for } \quad
A \in {\cal L}^{(k)}. $$
% To see this is a representation, consider two
% elementary tensors $L=A_1\otimes A_2\otimes\cdots\otimes A_k\otimes P_{k+1}\otimes
% P_{k+2}\otimes\cdots$ and
% $M=B_1\otimes B_2\otimes\cdots\otimes B_m\otimes P_{m+1}\otimes
% P_{m+2}\otimes\cdots$ in ${\cal L}$ and assume that $m>k$. Then
% $$ \pi_0(L)=\pi_m(A_1\otimes\cdots\otimes A_k\otimes P_{k+1}\otimes\cdots\otimes P_m)
% \,\P_m\quad\hbox{and}\quad
% \pi_0(M)=\pi_m(B_1\otimes\cdots\otimes B_m)\,\P_m, $$
% thus
% $$ \pi_0(L)\cdot\pi_0(M)=\pi_m(A_1B_1\otimes\cdots\otimes A_kB_k\otimes P_{k+1}B_{k+1}
% \otimes\cdots\otimes P_mB_m)\,\P_m =\pi_0(LM), $$
% using the fact that $\pi_m$ is a representation, and that it commutes with $\P_m$.
% By a similar argument we also get linearity, and the involutive property is trivial.
% To get norm continuity, observe that
% $$\big\|\pi_0(L)\big\|=\big\|\pi_k(A_1\otimes\cdots\otimes A_k)\,\P_{k+1}\big\|
% \leq\big\|A_1\otimes\cdots\otimes A_k\big\|=\|L\|\,.$$
% Thus $\pi_0$ is indeed a representation of ${\cal L}$.
That it is non-degenerate
follows from the fact that each $\pi_k$ is non-degenerate, and that
\break $\slim\limits_{k\to\infty}\P_k~=~\un$.
To see that $\widetilde\pi_0\rest{\cal A}
 =\pi$, recall that  $\pi_k$ is the representation
 obtained from from $\pi\rest{\cal A}^{(k)}$,
using the host algebra property of ${\cal L}^{(k)}\,.$
Let $B\in\cl A.^{(k)}$, then
for $A \in {\cal L}^{(k)}$ we have
$$ \wt\pi_0(B)
\,\pi_0(\Psi_k(A)) =
 \pi_0(B\cdot\Psi_k(A)) = \pi_k(B\cdot A)\P_{k+1}
=\pi(B)\,\pi_k(A)\P_{k+1}=\pi(B)\,\pi_0(\Psi_k(A))$$
from which it follows that $\widetilde\pi_0\rest{\cal A}
 =\pi\,.$
\qed

\def\bn{{\bf n}}

\def\bq{{\bf q}}

Thus for every family of projections $P_k\in{\cal L}_k$ we get a host algebra.
Now recall that ${\cal L}_k \cong {\cal K}(\ell^2(\N))$, and that
 there is a (countable) approximate
identity $(E_n)_{n \in \N}$ in ${\cal K}(\ell^2(\N))$ consisting of a strictly increasing
sequence of projections $E_n$ with $\dim (E_n\ell^2(\N)) = n$.
For each $k$, choose such an approximate identity $(E_n^{(k)})\subset{\cal L}_k$,
then for each sequence ${\bf n}=(n_1,\,n_2,\ldots)\in\N^\infty := \N^\N$,
we have a sequence of projections $\big(E_{n_1}^{(1)},\,E_{n_2}^{(2)},\ldots\big)$
from which we can construct an infinite tensor product as above, and we will
denote it by ${\cal L}[{\bf n}]$.
For the elementary tensors, we streamline the notation to:
$$A_1\otimes\cdots\otimes A_k\otimes E[{\bf n}]_{k+1}:=
A_1\otimes\cdots\otimes A_k\otimes E_{n_{k+1}}^{(k+1)}\otimes E_{n_{k+2}}^{(k+2)}
\otimes\cdots\in {\cal L}[{\bf n}]$$
  where $A_i\in{\cal L}_i$, and their closed span
is the simple $C^*$-algebra ${\cal L}[{\bf n}]$.

Next we want to define componentwise multiplication between different
C*-algebras ${\cal L}[{\bf n}]$ and ${\cal L}[{\bf m}]$.
This can of course be done in the algebraic infinite tensor product of the algebras ${\cal L}_k$,
(cf.~[Bo74, p470]) using suitable closures of subalgebras, but
it is faster to proceed as follows.
Note that  for componentwise multiplication, the sequences give:
$$\big(E_{n_1}^{(1)},\,E_{n_2}^{(2)},\ldots\big)\cdot\big(E_{m_1}^{(1)},\,E_{m_2}^{(2)},\ldots\big)
=\big(E_{p_1}^{(1)},\,E_{p_2}^{(2)},\ldots\big)$$
where $p_j:={\rm min}(n_j,\,m_j)$, i.e., multiplication reduces the entries, and hence
the sequence ${\big(E_1^{(1)},\,E_1^{(2)},\,E_1^{(3)}\ldots\big)}$ is invariant
under such multiplication. So we define an embedding
 ${\cal L}[{\bf n}]\subseteq M({\cal L}[{\bf 1}])$ for all ${\bf n}$,
where ${\bf 1}:=(1,\,1,\ldots)$ by
$$\eqalign{\big(A_1\otimes\cdots&\otimes A_k\otimes E[{\bf n}]_{k+1}\big)\cdot
\big(B_1\otimes\cdots\otimes B_n\otimes E[{\bf 1}]_{n+1}\big)  \cr
:=&\cases{A_1B_1\otimes\cdots\otimes A_nB_n\otimes A_{n+1}E_1^{(n+1)}\cdots\otimes A_kE_1^{(k)}\otimes E[{\bf 1}]_{k+1}
& if $n\leq k$ \cr
A_1B_1\otimes\cdots\otimes A_kB_k\otimes E_{n_{k+1}}^{(k+1)} B_{k+1}\cdots
\otimes E_{n_n}^{(n)} B_{n}\otimes E[{\bf 1}]_{n+1}
& if $n\geq k$ \cr}\cr}$$
for the left action, and similar for the right action on ${\cal L}[{\bf 1}]\,.$
To see that this is an embedding as claimed, choose a faithful representation $\pi_i$ of each
${\cal L}_i\cong{\cal K}({\cal H})$ on a Hilbert space ${\cal H}_i$ and
let $\psi_n$ be a unit vector in $E_1^{(n)}{\cal H}_n\,.$
Construct the infinite
tensor product Hilbert space $\bigotimes\limits_{n=1}^\infty{\cal H}_n$
w.r.t.\ the sequence ${(\psi_1,\,\psi_2,\ldots)}$, and note
that for each ${\cal L}[{\bf n}]$,
the tensor representation $\bigotimes\limits_{n=1}^\infty\pi_n$
on  $\bigotimes\limits_{n=1}^\infty{\cal H}_n$
is faithful (since it is faithful on the C*-algebras of which they are inductive limits).
Then it is obvious  that the given multiplication above is concretely realised
on this Hilbert space, and by faithfulness of the representations we realise the
embeddings ${\cal L}[{\bf n}]\subseteq M({\cal L}[{\bf 1}])$ for all ${\bf n}$.
Then
$${\cal L}[{\bf n}]\cdot{\cal L}[{\bf m}]
\subseteq{\cal L}[{\bf p}], \leqno(4.3) $$
where $p_j:={\rm min}(n_j,\,m_j)$, and in fact
$${\cal L}[{\bf n}]\subset M({\cal L}[{\bf p}])\supset
{\cal L}[{\bf m}].\leqno(4.4) $$

Using the embedding
${\cal L}[{\bf n}]\subseteq M({\cal L}[{\bf 1}])$ for all ${\bf n}$,
we define the
$C^*$-algebra in $M({\cal L}[{\bf 1}])$ generated by all ${\cal L}[{\bf n}]$,
and denote it by ${\cal L}[E]$. By (4.3), this is just the closed span
of all ${\cal L}[{\bf n}]$ and hence the closure of
the dense *-subalgebra $\cl L._0\subset{\cal L}[E]$, where
$$ \cl L._0:= \sum_{{\bf n}\in\N^\infty} \cl L.[{\bf n}]_0 \quad\hbox{(finite sums) and}\quad
\cl L.[{\bf n}]_0:=\bigcup_{k \in \N} {\cal L}^{(k)} \otimes E[\bn]_{k+1}. $$

We still have ${\cal A}\subset M({\cal L}[E])\supset{\cal L}^{(n)}$ for each
$n \in \N$.
Note that if two sequences ${\bf n}$ and ${\bf m}$ differ only  in a finite number of entries,
then ${\cal L}[{\bf n}]={\cal L}[{\bf m}]$, and hence we actually have
that the correct index set for the algebras ${\cal L}[{\bf n}]$ is
not the sequences $\N^\infty$, but the set of equivalence classes
$\N^\infty\big/\mathord{\sim}$ where ${\bf n}\sim{\bf m}$ if they
differ only in finitely many entries.
Some of the structures of $\N^\infty$ will factor through to
$\N^\infty\big/\mathord{\sim},$ e.g. we have a partial ordering
of equivalence classes defined by $[\b n.]\geq [\b m.]$ if
for any representatives $\b n.$ and $\b m.$ resp., we have that
there is an $N$ (depending on the representatives)
such that $n_k\geq m_k$ for all $k>N\,.$
In particular, we note that
products reduce sequences, i.e., we have $\cl L.[\b n.]\cdot\cl L.[\b p.]\subseteq
 \cl L.[\b q.]$ for $q_i={\rm min}(n_i,p_i),$ so $[\b n.]\geq[\b q.]\leq[\b p.]$.

Let $\varphi:\N^\infty\big/\mathord{\sim}
\to\N^\infty$ be a section of the factor map.
Then
 ${\cal L}[E]$ is the $C^*$-algebra generated in $M({\cal L}[{\bf 1}])$ by
${\big\{{\cal L}[\varphi(\gamma)]\,\big|\,\gamma\in\N^\infty\big/\mathord{\sim}
\big\}}$, and it is the closure of the span of the elementary tensors
in this generating set.

Below we will prove that ${\cal L}[E]$ is a full host algebra for
$(S,\sigma)$, and so it is of some interest to explore its algebraic structure.
{}From the reducing property of products, we already know that ${\cal L}[E]$
has the ideal ${\cal L}[{\bf 1}]$ (we will show that it is proper), hence
that it is not simple. However, it has in fact infinitely many proper ideals
and each of the generating algebras $\cl L.[\b n.]$ is contained in such an ideal:
\Proposition III.5.
For the $C^*$-algebra ${\cal L}[E]$, we have the following:
\item{(i)} ${\cal L}[E]$ is nonseparable,
\item{(ii)} Define $\cl I.[\b n._1,\ldots,\b n._k]$ to be the closed span of
$$\big\{ {\cal L}[{\bf q}]_0 \,\mid\, [\b q.]\leq[\b n._\ell]\;\hbox{for some}\;\;\ell=1,\ldots,k\big\}\,.$$
Let $[\b p.]>[\b n._\ell]$ strictly for all $\ell\in\{1,\ldots,\,k\}\,,$ then
$\cl L.[\b p.]\cap\cl I.[\b n._1,\ldots,\b n._k]=\{0\}\,.$
\item{(iii)} $\cl I.[\b n._1,\ldots,\b n._k]$ is a proper closed
two sided ideal of  ${\cal L}[E]\,.$
\item{(iv)} Define $\cl L.[\b n._1,\ldots,\b n._k]:= C^*\left(
 \cl L.[\b n._1]\cup\cdots\cup\cl L.[\b n._k]\right)\,.$ Then
 $\cl L.[\b n._1,\ldots,\b n._k]\subset\cl I.[\b n._1,\ldots,\b n._k]$
and
$$C^*\left(\cl L.[\b n._1,\ldots,\b n._k]
 \cdot\cl L.[\b n._{k+1}]\right)
\subseteq \cl L.[\b q._1,\ldots,\b q._k],\quad\hbox{where}\quad
 (\b q._{j})_\ell={{\rm min}\big((\b n._j)_\ell,\,
 (\b n._{k+1})_\ell\big)}\,.$$ % \leqno{(4.5)}$$

\Proof. (i) $\cl L.[E]\supset Q:=\big\{E[\b n.]_1:=E_{n_1}^{(1)}\otimes E_{n_2}^{(2)}
\otimes\cdots\,\mid\,\,\b n.\in\N^\infty\big\}\,.$
If $\b n.\not=\b p.$, there is some $k$ for which $E_{n_k}^{(k)}\not=E_{p_k}^{(k)}$
and as the approximate identity is linearly increasing, one of these must be
larger than the other, so take $E_{n_k}^{(k)}>E_{p_k}^{(k)}$ strictly.
Group the remaining parts of the tensor product together, i.e., write
$$ E[\b n.]_1=E_{n_k}^{(k)}\otimes A\quad\hbox{and}\quad
E[\b p.]_1=E_{p_k}^{(k)}\otimes B\,,   $$
where $A$ and $B$ are projections, then choose a product representation
$\pi=\pi_1\otimes\pi_2$ in which $\pi_1$ is faithful on $\cl L._k$ and
$\pi_2$ is faithful on the $C^*$-algebra
generated by $A$ and $B\,.$ Thus there is a unit vector $\psi\in\cl H._1$
such that $\big\|\pi_1(E_{n_k}^{(k)})\psi\big\|=1$ and $\pi_1(E_{p_k}^{(k)})\psi=0\,.$
For any unit vector $\phi\in\cl H._2$ we get
$$\leqalignno{\big\|E[\b n.]_1-E[\b p.]_1\big\| &\geq
\left\|(\pi_1\otimes\pi_2)\big(E_{n_k}^{(k)}\otimes A
-E_{p_k}^{(k)}\otimes B\big)(\psi\otimes\phi)\right\| \cr
&=\big\|\pi_1(E_{n_k}^{(k)})\psi\otimes\pi_2(A)\phi\big\|
=\big\|\pi_1(E_{n_k}^{(k)})\psi\big\|\cdot\big\|\pi_2(A)\phi\big\|
=\big\|\pi_2(A)\phi\big\| \cr}$$
and by letting $\phi$ range over the unit ball we get that
$\big\|E[\b n.]_1-E[\b p.]_1\big\|\geq\|A\|= 1\,.$
Thus, since $Q$ is uncountable and its elements far apart, $\cl L.[E]$ cannot be separable.

(ii) Here we adapt the argument in (i) as follows.
It suffices to show that for $\bq_1,\ldots, \bq_d$ with $\bq_i \leq \bn_j$ for some $j$,
the norm distance between
$\sum\limits_{i=1}^d\cl L.[\b q._i]_0$
and any $C\in \cl L.[\b p.]_0$ is always $\geq \|C\|\,.$
Let $C\in \cl L.[\b p.]_0$ be nonzero and consider a sum
$\sum\limits_{i=1}^d C_i$ with $C_i\in\cl L.[\bq_i]_0$
and $[\b p.]>[\b n._j]$ for all $j$, which implies
$[\b p.]>[\bq_i]$ for all $i$.
Choose an $M>0$ large enough so that all $C$ and $C_i$
can be expressed in
the form:
$$C_i=C_i^{(0)}\otimes E[\b n._i]_M\,,\quad\hbox{for}\quad
C_i^{(0)}\in\cl L.^{(M-1)}\,. $$
Then by $[\b p.]>[\bq_i]$
there is an entry of the tensor products, say for $j>M$,
which consist only of elements of the approximate identity
$(E_n^{(j)})_{n=1}^\infty\subset\cl L._j$ and
for which $B>B_i$ for all $i,$
where $B$ (resp. $B_i$) is the $j^{\rm th}$ entry of $C$  (resp. $C_i$).
Denote the remaining parts of the tensor products by $A$ (resp. $A_i$), i.e.,
$$  C=A\otimes B,\quad C_i=A_i\otimes B_i, \quad\hbox{where}\quad B>B_i\;\forall\,i$$
and $B,\; B_i$ consist of commuting projections. Then
$\|C-\sum\limits_{i=1}^d C_i\|=\big\|A\otimes B-\sum\limits_{i=1}^d
(A_i\otimes B_i)\big\|\,.$
% Recall that for the spatial tensor norm which we have here,
% $$\|A\otimes B\|\geq {\sup\big\{\|\pi(A\otimes B)\psi\|\,\mid\,\psi\in\cl H._\pi,\;
% \|\psi\|\leq 1,\;\pi\;\hbox{a product representation}\big\}}\,.$$
Choose a product representation $\pi=\pi_1\otimes\pi_2$ such that
$\pi_1$ is faithful on $\cl L.[\b p.]$ and $\pi_2$ is faithful on
the $C^*$-algebra generated by $(E_n^{(j)})_{n=1}^\infty\subset\cl L._j$.
Thus there is a unit vector $\phi\in\cl H._{\pi_2}$ such that $\|\pi_2(B)\phi\|=1$
and $\pi_2(B_i)\phi=0$ for all $i$ % (e.g. a unit vector in
% $\pi_2(B)\big(\mathop{\bigcup}\limits_{i=1}^n\pi_2(B_i)\cl H._{\pi_2}\big)^\perp$
(which exists because $B>B_i$ for all $i$). Then we have for any  unit vector
$\psi\in\cl H._{\pi_1}$ that
$$\eqalignno{\|C-\sum_{i=1}^d C_i\|&\geq\Big\|(\pi_1\otimes\pi_2)
\Big(A\otimes B-\sum_{i=1}^d A_i\otimes B_i\Big)(\psi\otimes\phi)\Big\|\cr
&=\|\pi_1(A)\psi\otimes\pi_2(B)\phi\|
=\|\pi_1(A)\psi\|\cdot\|\pi_2(B)\phi\|
=\|\pi_1(A)\psi\|\cr}$$
and by letting $\psi$ range over the unit ball of $\cl  H._{\pi_1}$,
we find that $\|C-\sum_{i=1}^d C_i\|\geq\|A\|=\|C\|$ since $\|B\|=1\,.$
This establishes the claim.

(iii) It is obvious from  the reduction property
$\cl L.[\b n.]\cdot\cl L.[\b p.]\subseteq
 \cl L.[\b q.]$ for $q_j={\rm min}(n_j,p_j),$ that ${\cl I.[\b n._1,\ldots,\b n._k]}$
is a two--sided ideal (hence a *--algebra). To see that it is proper, note
that $[\b p.]>[\b n._i]$ strictly for all $i$ where
$p_j={{\rm max}((\b n._1)_j,\ldots,(\b n._k)_j)+1}\,.$
Thus, by (ii) we see that
$\cl L.[\b p.]\cap\cl I.[\b n._1,\ldots,\b n._k]=\{0\}$ and hence that
${\cl I.[\b n._1,\ldots,\b n._k]}$ is proper.

(iv) ${\cl L.[\b n._1,\ldots,\b n._k]}\subset{\cl I.[\b n._1,\ldots,\b n._k]}$
because ${\cl I.[\b n._1,\ldots,\b n._k]}$ is a $C^*$-algebra which contains all
the generating elements ${\cl L.[\b n._i]}$ of ${\cl L.[\b n._1,\ldots,\b n._k]}\,.$
Next we need to prove that
$$C^*\left(\cl L.[\b n._1,\ldots,\b n._k]
 \cdot\cl L.[\b n._{k+1}]\right)
\subseteq \cl L.[\b q._1,\ldots,\b q._k],\quad\hbox{where}\quad
 (\b q._{j})_\ell={{\rm min}\big((\b n._j)_\ell,\,
 (\b n._{k+1})_\ell\big)}\,.$$
 By definition, $C^*\left(\cl L.[\b n._1,\ldots,\b n._k]
 \cdot\cl L.[\b n._{k+1}]\right)$ is the closed linear span of monomials
 $\prod\limits_{i=1}^NL_i$, where $L_i$ can be either of the form
 $A_iB_i$ or $B_iA_i$, where $A_i\in\cl L.[\b n._1,\ldots,\b n._k]$
 and $B_i\in\cl L.[\b n._{k+1}]\,.$ So it suffices to show that
 $$AB\in\cl L.[\b q._1,\ldots,\b q._k]\quad\hbox{for }\quad
 A\in\cl L.[\b n._1,\ldots,\b n._k]\quad\hbox{and}\quad
 B\in\cl L.[\b n._{k+1}]$$
 (since then $BA\in\cl L.[\b q._1,\ldots,\b q._k]$ by involution).
 Since $\cl L.[\b n.]_0$ is dense in $\cl L.[\b n.]$, it suffices to prove this for
 $A=A_1A_2\ldots A_p$ where $A_i=C_i\otimes E[\b n._{k_i}]_{r_i+1}$ and
 $C_i\in\cl L.^{(r_i)}\,,$ $k_i\in\{1,\ldots,k\}\,,$ and
 $B=D\otimes E[\b n._{k+1}]_{r+1}$, where $D\in\cl L.^{(r)}\,.$ Now
 $$ A_pB=F\otimes E[\b q._{k_p}]_{s+1}\in\cl L.[\b q._{k_p}]$$
 for some $F\in\cl L.^{(s)}\,,$ $s\geq{\rm max}(r_p,r)\,.$ Then
 $$ A_{p-1}A_pB=\big(C_{p-1}\otimes E[\b n._{k_{p-1}}]_{r_{p-1}+1}\big)
 \big(F\otimes E[\b q._{k_p}]_{s+1}\big)=G\otimes E[\b m.]_{t+1}, $$
where  $t\geq{\rm max}(r_{p-1},s)$ and
$$ \eqalign{  m_i&={\rm min}\big((\b n._{k_{p-1}})_i,\,(\b q._{k_p})_i)\big)
 ={\rm min}\Big((\b n._{k_{p-1}})_i,\,{\rm min}\big((\b n._{k_{p}})_i,
 \,(\b n._{k+1})_i\big)\Big) \cr
&={\rm min}\Big({\rm min}\big((\b n._{k_{p-1}})_i,\,(\b n._{k+1})_i\big),
 \,{\rm min}\big((\b n._{k_{p}})_i,\,(\b n._{k+1})_i\big)\Big)
 ={\rm min}\big((\b q._{k_{p-1}})_i,\,(\b q._{k_p})_i\big)\cr}$$
and so we have in fact that
$$ A_{p-1}A_pB=\big(\wt{C}\otimes E[\b q._{k_{p-1}}]_{t+1}\big)
 \big(\wt{F}\otimes E[\b q._{k_{p}}]_{t+1}\big)
\in\cl L.[\b q._{k_{p-1}}]\cdot\cl L.[\b q._{k_p}]$$
where $\wt{C}, \wt{F}\in\cl L.^{(t)}.$
Hence $A_{p-1}A_pB\in\cl L.[\b q._{k_{p-1}},\,\b q._{k_p}]\,.$
We continue the process to get $AB=A_1A_2\ldots A_pB\in\cl L.[\b q._1,\ldots,\b q._k]\,.$
\qed
\noindent
For each strictly increasing sequence $\big([\b n._1],\,[\b n._2],\ldots\big)\subset
\N^\infty/\mathord{\sim}$ we get from part (ii) a strictly increasing chain of proper ideals
$\cl J._k:=\cl I.[\b n._1,\ldots,\b n._k]\,.$

Now we want to prove our main theorem in this section.
\Theorem III.6. The monomorphism
$\eta:S_\sigma\to U(M({\cal L}[E]))$ from above, defined by
$$\eta((s,t)):= t\delta_s\in {\cal A}\subset M({\cal L}[E]),$$ is
continuous with respect to the strict topology on $M({\cal L}[E])$
and
${\cal L}[E]$ is a host algebra, i.e., the map
$$\eta^*:{\rm Rep}({\cal L}[E],\cl H.)\to {\rm Rep}((S,\sigma),\cl H.)$$
is injective. The range of $\eta^*$ is exactly $\cal R(\cl H.)$.

\Proof. First we show that
$\eta$ is continuous with respect to the strict topology on
$M({\cal L}[E])$. This implies that for each $\pi\in{\rm Rep}({\cal L}[E],\cl H.)$
the representation $\widetilde\pi\in{\rm Rep}({\cal A},\cl H.)$ is regular, hence
$$\eta^*\big({\rm Rep}({\cal L}[E],\cl H.)\big)\subseteq{\cal R}(\cl H.).$$
Since $\im(\eta)$ is bounded, it suffices to show that the set
$$\big\{L \in {\cal L}[E]\,\mid\,g\mapsto
\eta(g)L\quad\hbox{is norm continuous in}\quad g\in
S_\sigma\big\}$$
spans a dense subspace of ${\cal L}[E]$. This reduces the assertion
to the corresponding result for the action of $S_\sigma$ on ${\cal L}[\bn]$ for
each $\bn$, which follows from the continuity of the corresponding map
$S_\sigma \to M({\cal L}[\bn])$ (Proposition~III.4).

%For each ${\cal L}[{\bf n}]\subset{\cal L}[E]$, denote its essential subspace
%by ${\cal H}_{\bf n}:=\pi\big({\cal L}[{\bf n}]\big){\cal H}$. Then we can apply
%Proposition~IV.4 to the restricted representation
%$\pi_{\bf n}(L):=\pi(L)\rest{\cal H}_{\bf n}$, where $L\in{\cal L}[{\bf n}]$
%to conclude that $\widetilde\pi_{\bf n}\rest{\cal A}=
%\left(\widetilde\pi\rest{\cal A}\right)\rest{\cal H}_{\bf n}$
%is regular. But the ${\cal L}[{\bf n}]$ generate  ${\cal L}[E]$ as the ${\bf n}$
%varies, hence ${\cal H}_{\bf n}$ span ${\cal H}$ as ${\bf n}$ varies.
%Thus $\widetilde\pi\in{\rm Rep}({\cal A},\cl H.)$ is regular.

To prove that $\eta^*$ is injective
we show that $\cal A$ separates ${\rm Rep}({\cal L}[E],\cl H.)$ for all $\cl H.$.
Let $\pi\in{\rm Rep}({\cal L}[E],\cl H.)$, then by
Proposition~III.4 we know that the values which $\widetilde\pi(\cal A)$
takes on ${\cal H}_{\bf n}$ uniquely determine the values of $\pi\big({\cal L}[{\bf n}])$
on its essential subspace ${\cal H}_{\bf n}$, hence on all $\cal H$, as $\pi\big({\cal L}[{\bf n}]\big)$
is zero on the orthogonal complement of ${\cal H}_{\bf n}$.
This holds for all ${\bf n}$, hence $\widetilde\pi(\cal A)$ uniquely determines the values of
$\pi$ on ${\cal L}[E]$, i.e., $\eta^*$ is injective.

 It remains to prove that
 $\eta^*\big({\rm Rep}({\cal L},\cl H.)\big)={\cal R}(\cl H.)$. Start from a
 $\pi\in{\rm Rep}({\cal A},\cl H.)$ which is regular. Then we have to show how to obtain
 a $\pi_0\in{\rm Rep}\,{\cal L}[E]$ such that $\widetilde\pi_0\rest{\cal A}
 =\pi$.
 Observe that $\pi$ is regular on all ${\cal A}^{(n)}$, hence there are unique
 $\pi_n\in{\rm Rep}({\cal L}^{(n)},\cl H.)$ which extend (on $\cl H.$) to
coincide with
 $\pi\rest{\cal A}^{(n)}$ by the host algebra property of ${\cal L}^{(n)}$.
For each $\bf n$ define the projections
$$ \E^{\bf n}_{k}:= \slim_{m\to\infty}\pi(E_{n_k}^{(k)})\cdots\pi(E_{n_m}^{(m)})
\quad\hbox{and}\quad \E^{\bf n}:=\slim_{k\to\infty}\E^{\bf n}_{k}\,.$$
Now each $\pi_n({\cal L}^{(n)})$ commutes with the projections $\E^{\bf n}_{k}$
for $k>n$, and in particular
preserves the space ${\cal H}^{\bf n}:=\E^{\bf n}{\cal H}$,
and hence so does $\pi({\cal A}^{(n)})$.
Then by Proposition~III.4 we know that we can define a (non-degenerate) representation
$\pi_0^{\bf n}:{\cal L}[{\bf n}]\to{\cal B}({\cal H}^{\bf n})$ by
$$\pi_0^{\bf n}(L)
=\pi_k(A_1\otimes\cdots\otimes A_k)\,\E^{\bf n}_{k+1}$$
for $L=A_1\otimes\cdots\otimes A_k\otimes E_{n_{k+1}}^{(k+1)}\otimes E_{n_{k+2}}^{(k+2)}
\otimes\cdots\in {\cal L}[{\bf n}]$ such that $\widetilde\pi_0^{\bf n}
\rest{\cal A}$
is $\pi({\cal A})$, restricted to ${\cal H}^{\bf n}$.
We extend $\pi_0^{\bf n}$ to all of ${\cal H}$, by putting it to zero on
the orthogonal complement of ${\cal H}^{\bf n}$. Note that
$$ {\bf n} \leq {\bf m} \quad \Rarrow \quad {\cal H}^{\bf n} \subeq {\cal H}^{\bf m}. $$

We now argue that these representations $\pi_0^{\bf n}$ combine into a single
representation of ${\cal L}[E]$.
First, we want to extend by linearity
the maps $\pi_0^{\bf n}:{\cal L}[{\bf n}]
\to{\cal B}({\cal H})$ to define a linear map $\pi_0$ from
the dense $*$-subalgebra $\cl L._0\subset{\cal L}[E]$ to ${\cal B}({\cal H})$, where
we recall that
$ \cl L._0:= \sum\limits_{{\bf n}\in\N^\infty} \cl L.[{\bf n}]_0$
(finite sums).
%  \quad\hbox{where}\quad
% \cl L.[{\bf n}]_0:=\bigcup_{k \in \N} {\cal L}^{(k)} \otimes E[\bn]_{k+1}. $$

This linear extension $\pi_0$ is possible if the sum of the spaces $\cl L.[{\bf n}]_0$ is
direct for different ${\bf n}\in\varphi\big(\N^\infty\big/\mathord{\sim}\big)$,
i.e., if $0=\sum\limits_{k=1}^m B_k$ for $B_k\in
\cl L.[{\bf n}_k]_0,$ where ${\bf n}_k\not\sim\b n._\ell$ if $k\not=\ell$
implies that $B_k=0$ for all $k\,.$
Let us prove this implication, so assume  $0=\sum\limits_{k=1}^m B_k$
as above.
Choose an $M>0$ large enough so that for all $k,$ the $B_k$ can be expressed in
the form
$B_k=A_k\otimes E[\b n._k]_M$ for
$A_k\in\cl L.^{(M-1)}\,, $
define the projections $P_k:=\un\otimes\cdots\otimes\un\otimes
 E[\b 1.]_k$ (there are $k-1$ factors of $\un$), and note that $P_\ell$ commutes
 with all $B_k$ for $\ell\geq M.$ In fact, for $B_k$ as above, we have
 (simplifying notation to $\b n._k=\b n.$):
$$B_kP_\ell=A_k\otimes
E_{n_M}^{(M)}\otimes\cdots\otimes E_{n_{\ell-1}}^{(\ell-1)}
\otimes  E[\b 1.]_\ell \in \cl L.^{(\ell-1)}\otimes E[\b 1.]_\ell
$$
and so multiplication by $P_\ell$ for $\ell\geq M$ maps the
$B_k$ to elementary tensors of the form
$A_k\otimes
E_{n_M}^{(M)}\otimes\cdots\otimes E_{n_{\ell-1}}^{(\ell-1)}$ in $\cl L.^{(\ell-1)}$
(after identifying $\cl L.^{(\ell-1)}\otimes E[\b 1.]_\ell$
with $\cl L.^{(\ell-1)}$).
Now a set of elementary tensors (in a finite tensor product) will be linearly
independent if the entries in a fixed slot are linearly independent
so it suffices to find $\ell>M$ such that
the pieces ${E_{n_M}^{(M)}\otimes\cdots\otimes E_{n_{\ell-1}}^{(\ell-1)}}$
are linearly independent for $\b n.\in N:={\big\{\b n._k\,\mid\,
k=1,\ldots,\,m\big\}}.$
Since the approximate identities
$(E_n^{(k)})_{n=1}^\infty\subset\cl L._k$ consist of strictly increasing
projections, their terms are linearly independent from which it follows
that tensor products of these with distinct entries are linearly independent.
Thus we only have to identify an $\ell$ large enough so that the portions of the
sequences $\b n._k$ between the entries $M$ and $\ell$ can distinguish
all the sequences in $N\,,$ and this is always possible since
the $\b n._k$ are representatives of distinct
equivalence classes in $\N^\infty/\mathord{\sim}$.
Thus $\{B_1P_\ell,\ldots,\,B_mP_\ell\}$ is linearly independent for this $\ell,$
so $0=\sum\limits_{k=1}^m B_kP_\ell$ implies that all $B_k=0$.
We conclude that the linear extension $\pi_0$ exists.

 That $\pi_0$ respects involution is clear. To see that it is
a homomorphism, consider the
elementary tensors $$L=A_1\otimes A_2\otimes\cdots\otimes A_k\otimes E[{\bf n}]_{k+1}
\in{\cal L}[{\bf n}]\quad\hbox{and}\quad
M={B_1\otimes B_2\otimes\cdots\otimes B_m\otimes E[{\bf p}]_{m+1}\in
{\cal L}[{\bf p}]}$$
where $m>k$ and ${\bf n} \not\sim {\bf p}\in \N^\infty$. Then
$$\eqalignno{\pi_0(L)\,\pi_0(M)&=
\pi_k(A_1\otimes\cdots\otimes A_k)\,\E^{\bf n}_{k+1}
\pi_m(B_1\otimes\cdots\otimes B_m)\,\E^{\bf p}_{m+1}\cr
&=\pi_m\big(A_1\otimes\cdots\otimes A_k\otimes E_{n_{k+1}}^{(k+1)}\otimes\cdots
\otimes E_{n_{m}}^{(m)}\big)\,\E^{\bf n}_{m+1}
\pi_m(B_1\otimes\cdots\otimes B_m)\,\E^{\bf p}_{m+1}\cr
&=\pi_m\big(A_1B_1\otimes\cdots\otimes  A_kB_k\otimes E_{n_{k+1}}^{(k+1)}B_{k+1}
\cdots\otimes E_{n_{m}}^{(m)} B_m\big)\,\E^{\bf n}_{m+1}\,\E^{\bf p}_{m+1}\,.}$$
Now recall that the operator product is jointly continuous on bounded sets
in the strong operator topology, hence
$$\eqalignno{\E^{\bf n}_{k}\,\E^{\bf p}_{k}&=
\slim_{m\to\infty}\pi(E_{n_k}^{(k)})\cdots\pi(E_{n_m}^{(m)})
\cdot \slim_{r\to\infty}\pi(E_{p_k}^{(k)})\cdots\pi(E_{p_r}^{(r)})\cr
&=\slim_{m\to\infty}\pi(E_{n_k}^{(k)})\cdots\pi(E_{n_m}^{(m)})
\pi(E_{p_k}^{(k)})\cdots\pi(E_{p_m}^{(m)})\cr
&=\slim_{m\to\infty}\pi(E_{q_k}^{(k)})\cdots\pi(E_{q_m}^{(m)})
=\E^{\bf q}_{k}}$$
where $q_j:={\rm min}(n_j,\,p_j)$. Thus we get exactly that
$\pi_0(L)\,\pi_0(M)=\pi_0(LM)$.

We now verify that $\pi_0$
is bounded.
 For this, we first need to prove the following:
 \medskip\noindent
 \b Claim:. {\sl Recall that $\cl L.[\b n._1,\ldots,\b n._k]= C^*\left(
 \cl L.[\b n._1]\cup\cdots\cup\cl L.[\b n._k]\right)\,.$
 Then for each $k\geq 1$ and $k$-tuple ${(\b n._1,\ldots,\b n._k)}$
 such that ${\bf n}_k\not\sim\b n._\ell$ if $k\not=\ell$ the map
 $\pi_0$ on $\cl L._0\cap\cl L.[\b n._1,\ldots,\b n._k]$
 extends to a representation of the $C^*$-algebra
 $\cl L.[\b n._1,\ldots,\b n._k]\,.$}
\medskip\noindent
 \b Proof:.
 Note that the claim implies
 the compatibility of the representations, i.e.,
 on intersections $\cl L.[\b p._1,\ldots,\b p._\ell]\cap\cl L.[\b n._1,\ldots,\b n._k]\,,$
 the representations
 produced by $\pi_0$ on $\cl L.[\b n._1,\ldots,\b n._k]$
 and $\cl L.[\b p._1,\ldots,\b p._\ell]$ coincide.
 This is because $\pi_0$ is given as a consistent map on the dense space $\cl L._0\,.$\chop
We now prove the claim by induction on $k.$
 We already have
 by definition that $\pi_0$
 is the representation $\pi^{\b n.}$ on $\cl L.[\b n.]$ for each
 $\b n.,$ hence the claim is true for $k=1\,.$
 Assume the claim is true for all values of $k$ up to
 a fixed $k\geq 1,$ then we now prove it
 for $k+1\,.$ Observe that $\cl L.[\b n._1,\ldots,\b n._{k+1}]$
 contains the closed two--sided ideals
 $$\leqalignno{\cl J._1:=&C^*\left(\cl L.[\b n._1,\ldots,\b n._k]
 \cdot\cl L.[\b n._{k+1}]\right)\subset\cl J._2\cap\cl J._3, \cr
 \cl J._2:=& \cl J._1+\cl L.[\b n._1,\ldots,\b n._k]\quad\hbox{and}\quad
 \cl J._3:=\cl J._1+\cl L.[\b n._{k+1}]&\hbox{where}\cr}$$
 and that $\cl L.[\b n._1,\ldots,\b n._{k+1}]=\cl J._2 + \cl J._3$.
 We will prove below that $\cl J._1$ is proper (hence that the ideal structure above is nontrivial).
 Consider the factorization $\xi:\cl L.[\b n._1,\ldots,\b n._{k+1}]\to
 \cl L.[\b n._1,\ldots,\b n._{k+1}]\big/\cl J._1\,.$
 Then
 $$\xi\big(\cl L.[\b n._1,\ldots,\b n._{k+1}]\big)=\xi\big(\cl L.[\b n._1,\ldots,\b n._{k}]\big)
 +\xi\big(\cl L.[\b n._{k+1}]\big)$$
 and $\xi(\cl J._2)\cdot\xi(\cl J._3)=0\,.$
 If $\cl J._1$ is {\it not} proper, then
 $$\cl L.[\b n._{k+1}]\subset\cl J._1\supset\cl L.[\b n._1,\ldots,\b n._{k}]\,.$$
% Now products reduce sequences, i.e., we have $\cl L.[\b n.]\cdot\cl L.[\b p.]\subseteq
% \cl L.[\b q.]$ for $q_i={\rm min}(n_i,p_i),$ so $[\b n.]\geq[\b q.]\leq[\b p.]$.
By Proposition~III.5(iv), we have that
 $$\cl J._1\subset\cl L.[\b q._1,\ldots,\b q._k]\subset\cl I.[\b q._1,\ldots,\b q._k]
\quad \hbox{ for }\quad (\b q._{j})_\ell=  {{\rm min}\big((\b n._j)_\ell,\,
 (\b n._{k+1})_\ell\big)}\,,$$
% algebras $\cl L.[\b q.]$ for $[\b q.]\leq[\b q._\ell]$ for some $\ell=1,\ldots,k\,,$
and hence  $\cl L.[\b n._{k+1}]\subset\cl J._1\subset\cl I.[\b q._1,\ldots,\b q._k]\,.$
Thus, by Proposition~III.5(ii) we conclude that $[\b n._{k+1}]$ cannot be strictly greater
than all the $[\b q._i]\,,$ i.e., there is
one member of the set ${\{\b q._1,\ldots,\b q._k\}},$ say $\b q._j,$
which satisfies $[\b q._j]=[\b n._{k+1}],$ and so by definition of $\b q._j$, we have that
  eventually $(\b n._{k+1})_\ell={{\rm min}\big((\b n._j)_\ell,\,
 (\b n._{k+1})_\ell\big)},$ i.e., $[\b n._j]\geq[\b n._{k+1}]\,.$

 Likewise, the inclusion $\cl L.[\b n._1,\ldots,\b n._{k}]\subset
 {C^*\left(\cl L.[\b n._1,\ldots,\b n._k]\cdot\cl L.[\b n._{k+1}]\right)}=\cl J._1$
 implies that no $\b n._j,$ $j=1,\ldots,k,$ is reduced
 through multiplication by $\b n._{k+1},$ i.e.,
 eventually $(\b n._{j})_\ell={{\rm min}\big((\b n._j)_\ell,\,
 (\b n._{k+1})_\ell\big)}$ for all $j,$  i.e., $[\b n._j]\leq[\b n._{k+1}]\,.$
 So, together with the previous inequality, we see that there must be a $j\in\{
 1,\ldots,k\}$ such that $[\b n._j]=[\b n._{k+1}]\,.$ This contradicts the
 initial assumption that all $[\b n._\ell]$ are distinct, and so $\cl J._1$
 must be proper.

 Now consider $\pi_0$ on $\cl L._0\cap\cl L.[\b n._1,\ldots,\b n._{k+1}]\,.$
 By the induction assumption,  $\pi_0$ on
$$\cl L._0\cap\cl L.[\b n._1,\ldots,\b n._{k}]$$
 is the restriction of a representation on $\cl L.[\b n._1,\ldots,\b n._{k}]$,-
 we denote the projection onto its essential subspace by
 $\E[\b n._1,\ldots,\b n._{k}]\,.$ Note that $\E[\b n._{k+1}]$ commutes with
 $\E[\b n._1,\ldots,\b n._{k}]$ because it commutes with all the generating
 elements $\pi_0(L_i)=\pi^{\b n._i}(L_i)\,,$ $L_i\in\cl L.[\b n._i]$.
 Thus we have an orthogonal decomposition $\cl H.=\cl H._1\oplus\cl H._2
 \oplus\cl H._3\oplus\cl H._4$, where
 $$\eqalignno{\cl H._1:=&\E[\b n._1,\ldots,\b n._{k}]\E[\b n._{k+1}]\cl H.\,,\qquad
 \cl H._2:=\E[\b n._1,\ldots,\b n._{k}]\big(\un-\E[\b n._{k+1}]\big)\cl H.\cr
 \cl H._3:=&\E[\b n._{k+1}]\big(\un-\E[\b n._1,\ldots,\b n._{k}]\big)\cl H.\,,\qquad
 \cl H._4:=\big(\un-\E[\b n._{k+1}]\big)\big(\un-\E[\b n._1,\ldots,\b n._{k}]\big)\cl H.\cr}$$
 and $\pi_0$ preserves these subspaces.
 Now by Proposition~III.5(iv) and the induction assumption, $\pi_0$ extends from the
 $\cl L._0\cap\cl J._1$ to a representation on $\cl J._1,$ and as
$\cl J._1=C^*\left(\cl L.[\b n._1,\ldots,\b n._k]\cdot\cl L.[\b n._{k+1}]\right)$,
the essential projection for $\pi_0\restriction\cl J._1$ is
$\E[\b n._1,\ldots,\b n._{k}]\E[\b n._{k+1}]$, i.e., its essential subspace
is $\cl H._1\,.$ But since $\cl J._1$ is a closed two--sided ideal of
$\cl L.[\b n._1,\ldots,\b n._{k+1}]$, its non-degenerate representations
extend uniquely to $\cl L.[\b n._1,\ldots,\b n._{k+1}]\,.$
Thus on $\cl H._1,$ $\pi_0$ extends from $\cl L._0\cap\cl L.[\b n._1,\ldots,\b n._{k+1}]$
to a representation on $\cl L.[\b n._1,\ldots,\b n._{k+1}]\,.$

Next observe that on $\cl H._1^\perp=\cl H._2\oplus\cl H._3\oplus\cl H._4$ we have
$\{0\}=\pi_0(\cl J._1)$.
% If $\cl L.[\b n._1,\ldots,\b n._{k}]\subset\cl J._1$ then $\pi_0$ on
% $\cl H._2$ just becomes the zero representation. Otherwise,
We show that one can define a consistent representation of
$\xi\big(\cl L.[\b n._1,\ldots,\b n._{k+1}]\big)$ by $\rho(\xi(A))
:=\pi_0(A)\restriction\cl H._1^\perp,$ for $A\in\cl L.[\b n._{k+1}]
+ \cl L.[\b n._1,\ldots,\b n._{k}]$, using the structure of $\xi\big(\cl L.[\b n._1,\ldots,
\b n._{k+1}]\big)$ above. First observe that $\rho$ is well-defined on
$\xi\big(\cl L.[\b n._{k+1}]\big)$ and $\xi\big(\cl L.[\b n._1,\ldots,\b n._{k}]\big)$ separately,
because if $A_1-A_2\in\cl J._1$, then $\pi_0(A_1-A_2)\restriction\cl H._1^\perp=0\,.$
Next,  $\rho$ is well-defined on the set
$\xi\big(\cl L.[\b n._{k+1}]+\cl L.[\b n._1,\ldots,\b n._{k}]\big)$
by the induction assumption, and the
 consistency of the extensions of $\pi_0\,.$
To see that  $\rho$ is well-defined on the algebra
$\xi\big(\cl L.[\b n._1,\ldots,\b n._{k+1}]\big)=\xi\big(\cl L.[\b n._1,\ldots,\b n._{k}]\big)
 +\xi\big(\cl L.[\b n._{k+1}]\big),$
it suffices by the direct sum decomposition
to check it on $\cl H._2,$ $\cl H._3$ and $\cl H._4$ separately.
On $\cl H._2$, $\pi_0$ vanishes on $\cl L.[\b n._{k+1}]$,
so since $\xi\big(\cl L.[\b n._{k+1}]\big)$
is an ideal of $\xi\big(\cl L.[\b n._1,\ldots,\b n._{k+1}]\big)$
(and $\xi(\cl J._2)\cdot\xi(\cl J._3)=\{0\}$), it follows that we can extend
 $\rho(\xi(A))\restriction\cl H._2$ by linearity, i.e.,
 $\rho(\xi(A)+\xi(B))=\rho(\xi(A))$ for $A\in\cl L.[\b n._1,\ldots,\b n._{k}],$
 $B\in\cl L.[\b n._{k+1}]$
 to define a representation on $\xi\big(\cl L.[\b n._1,\ldots,\b n._{k+1}]\big)\,.$
Likewise, on
$\cl H._3,$ $\pi_0$ vanishes on $\cl L.[\b n._1,\ldots,\b n._{k}]$,
so we can show $\rho$ defines a representation of
$\xi\big(\cl L.[\b n._1,\ldots,\b n._{k+1}]\big)$
and on
$\cl H._4,$ $\rho$ is zero.
Then $\rho$ lifts to a representation of $\cl L.[\b n._1,\ldots,\b n._{k+1}]$
on $\cl H._1^\perp$
which coincides with $\pi_0$ on $\cl L._0\cap\cl L.[\b n._1,\ldots,\b n._{k+1}]\,.$
Taking the direct sum of this with the representation we obtained on $\cl H._1,$
produces a representation of $\cl L.[\b n._1,\ldots,\b n._{k+1}]$ on all $\cl H.$
which coincides with $\pi_0$ on $\cl L._0\cap\cl L.[\b n._1,\ldots,\b n._{k+1}]\,.$
Thus, we have proven the claim for $k+1$, which completes the induction.
 {}\hfill$\blacktriangledown$\break

That $\pi_0$ is bounded on $\cl L._0$ now follows immediately from the claim, because
any $A\in\cl L._0$ is of the form $A=\sum\limits_{k=1}^m B_k$
 for  $B_k\in
 \cl L.[{\bf n}_k]_0$, where ${\bf n}_k\not\sim\b n._\ell$ if $k\not=\ell\,.$
But this is an element of $\cl L.[\b n._1,\ldots,\b n._m]$ and by the claim
$\pi_0$ extends as a representation to it, hence $\|\pi_0(A)\|\leq\|A\|\,.$
We conclude that
$\pi_0$ is a bounded representation, hence extends to all of $\cl L.[E]\,.$
To see that $\pi_0$ is non-degenerate, recall that $\{E_n^{(k)}\}\subset\cl L._k$
is an approximate identity of increasing projections.
Thus we can find a sequence $\bf n$ such that
$\slim\limits_{m\to\infty}\pi(E_{n_m}^{(m)})=\un  $, and hence
$\E^{\bf n}=\un$ by $\pi(E_{n_m}^{(m)})\leq\E^{\bf n}\leq\un$
for all $m\,.$ Since the essential subspace of $\pi_0\rest\cl L.[{\bf n}]$
is $\E^{\bf n}\cl H.,$ it follows that $\pi_0$ is non-degenerate.
It then follows from Proposition~III.4 applied to $\cl L.[{\bf n}]$
that $\widetilde\pi_0\rest{\cal A}
 =\pi$.
\qed

% The following lemma shows that our seemingly special choice of $(S,B)$
% covers all countably dimensional symplectic vector spaces:
%
% \Lemma IV.7. In each countably dimensional symplectic vector space
% $(V,\Omega)$, there exists a basis $(p_n,q_n)_{n \in \N}$ with
% $$ \Omega(p_n,q_m) = \delta_{nm} \quad \hbox{ and } \quad
% \Omega(p_n,p_m) = \Omega(q_n,q_m) = 0 \quad \hbox{ for } \quad n,m \in \N. $$
% Then $Ip_n := -q_n$ and $Iq_n = p_n$ defines a complex structure on $V$
% for which
% $(v,w) := \Omega(Iv,w)$ is positive definite and extends to a
% positive definite hermitian form on $(V,I)$ with
% $$ \Omega(v,w) = -\Re(Iv,w) = \Im(v,w). $$
% % \qeddis
%
% \Proof. Let $(e_n)_{n \in \N}$ be a linear basis of $V$.
% We construct the basis elements $p_n, q_n$ inductively as follows.
% If $p_1,\ldots, p_k$ and $q_1, \ldots, q_k$ are already chosen, pick a minimal
% $m$ with $e_m \not\in \span\{p_1,\ldots, p_k, q_1,\ldots, q_k\}$ and put
% $$ p_{k+1} := e_m - \sum_{i=1}^k\big( B(e_m, q_i)p_i +B(p_i, e_m) q_i\big) $$
% to ensure that this element is $B$-orthogonal to all previous ones.
% Then pick $\ell$ minimal, such that $B(p_{k+1}, e_\ell) \not=0$, put
% $$ \tilde q_{k+1} := e_\ell - \sum_{i=1}^k\big( B(e_\ell, q_i)p_i +B(p_i, e_\ell) q_i\big) $$
% and pick $q_{k+1} \in \R \tilde q_{k+1}$ with
% $B(p_{k+1}, q_{k+1}) = 1$.
% This process can be repeated ad infinitum and produces the required
% bases of $V$ because for each $k$, the span of
% $p_1,\ldots, p_k, q_1,\ldots, q_k$ contains at least
% $e_1,\ldots, e_k$.
% \qed

Finally, we apply the structures above to produce a direct integral of
regular representations into irreducible regular representations.
First observe that given any representation $\pi\in
{\rm Rep}((S,\sigma),\cl H.)$, where $\cl H.$ is separable,
then as ${(E_n)_{n\in\N}}$ is an approximate identity
for $\cl K.(\ell^2(\N)),$ there is a sequence ${\bf n}$ such that
$\slim\limits_{k\to\infty}\slim\limits_{\ell\to\infty}\pi(E_{n_k}^{(k)})\cdots
\pi(E_{n_\ell}^{(\ell)})=\un,$ and thus by Proposition~III.4 there is
a unique $\pi_0\in{\rm Rep}(\cl L.[{\bf n}],\cl H.)$ such that $\eta^*\pi_0=\pi.$
Fix a choice of maximally commutative subalgebra $\cl C.\subset \pi_0(\cl L.[{\bf n}])'.$
Then, since $\cl L.[{\bf n}]$ is separable, there is an extremal decomposition
of $\pi_0$ (cf.~[BR03] Corollary~4.4.8), i.e., there is a standard measure space
${(Z,\mu)}$ with $\mu$ a positive bounded measure, a measurable family $z\to\cl H.(z)$
of Hilbert spaces, a measurable family $z\to\pi_z\in{\rm Rep}(\cl L.[{\bf n}],\cl H.(z))$
of representations which are almost all irreducible and a unitary
$U:\cl H.\to\int_Z^\oplus\cl H.(z)\,d\mu(z)$
such that $U\cl C.U^{-1}$ is the diagonizable operators, and
$$U\pi_0(A)U^{-1}=\int_Z^\oplus\pi_z(A)\,d\mu(z)\;\;\forall\,A\in\cl L.[{\bf n}].$$
Then for $\psi,\,\phi\in\int_Z^\oplus\cl H.(z)\,d\mu(z)$ with
decompositions $\psi=\int_Z^\oplus\psi_z\,d\mu(z)$
and $\phi=\int_Z^\oplus\phi_z\,d\mu(z)$,
we have for $s\in S$ and any countable approximate identity $(F_k)$ of $\cl L.[{\bf n}]$ that
$$
\leqalignno{&\big(\phi,\,U\pi(s)U^{-1}\psi\big)=\big(\phi,\,U\eta^*\pi_0(s)U^{-1}\psi\big)
=\lim_{k\to\infty}\big(\phi,\,U\pi_0(\delta_sF_k)U^{-1}\psi\big)\cr
&=\lim_{k\to\infty}\int_Z\big(\phi_z,\,\pi_z(\delta_sF_k)\,\psi_z\big)\,d\mu(z)
=\int_Z\lim_{k\to\infty}\big(\phi_z,\,\pi_z(\delta_sF_k)\,\psi_z\big)\,d\mu(z)\cr
&=\int_Z\big(\phi_z,\,\eta^*\pi_z(s)\,\psi_z\big)\,d\mu(z)
=\big(\phi,\,\int_Z^\oplus\eta^*\pi_z(s)\,d\mu(z)\psi\big),}
$$
where the usage of the Dominated Convergence Theorem in the second line is
justified by $\big|\big(\phi_z,\,\pi_z(\delta_sF_k)\,\psi_z\big)\big|\leq
\|\phi_z\|\|\psi_z\|$ as both of $z\to\phi_z$ and $z\to\psi_z$  are square integrable
w.r.t. $\mu.$ Hence
$$U\pi(s)U^{-1}=\int_Z^\oplus\eta^*\pi_z(s)\,d\mu(z)\;\;\forall\,s\in S.$$
Since $\eta^*$ preserves irreducibility, almost all $\eta^*\pi_z$ are irreducible,
and hence we obtain the promised decomposition.

\sectionheadline{Appendix. Host algebras and the strict topology}

\Lemma A.1. Let $X$ be a locally compact space.

{\rm(a)} On each bounded subset of
$M(C_0(X)) \cong C_b(X)$, the strict topology coincides with the topology of
compact convergence, i.e., the compact open topology.
This holds in particular for the subgroup $C(X,\T) \cong U(C_b(X))$.

{\rm(b)} A unital $*$-subalgebra $S \subeq C_b(X)$ is strictly dense if and only if
it separates the points of $X$.

\Proof. (a) ([Bl98, Ex.~12.1.1(b)]) Let
${\cal B} \subeq C_b(X)$ be a bounded subset with $\|f\| \leq C$ for each
$f \in {\cal B}$.
For each $\phi \in C_0(X)$ and $\eps > 0$ we now find a compact subset $K \subeq X$ with
$|\phi| \leq \eps$ outside $K$. For $f_i \to f$ in ${\cal B}$ with respect to the compact
open topology, we then have
$$ \|(f-f_i)\phi\|
\leq \|(f-f_i)\res_K\| \|\phi\|+ \eps\|f-f_i\|
\leq \eps \|\phi\| + 2\eps C $$
for sufficiently large $i$. Therefore the maps
${\cal B} \to C_0(X), f \mapsto f\phi$ are continuous if ${\cal B}$ carries the compact
open topology. This means that the strict topology on ${\cal B}$ is coarser than the
compact open topology.

If, conversely, $K \subeq X$ is a compact subset and $h \in C_0(X)$ with
$h\res_K = 1$, then
$$ \|(f - f_i)\res_K\| \leq \|(f-f_i)h\| $$
shows that the strict topology on $C_b(X)$ is finer than the compact open topology.
This proves~(a).

(b) If $S$ is strictly dense, then it obviously separates the points of
$X$ because the point evaluations are strictly continuous.

Suppose, conversely, that $S$ separates the points of $X$.
Replacing $S$ by its norm closure, we may w.l.o.g.\ assume that $S$ is
norm closed. Let $K \subeq X$ be compact. Since $S$ separates the points of $K$,
the Stone-Weierstra\3 Theorem implies that
$S\res_K = C(K)$. For any $f \in C_b(X)$ we therefore find some
$f_K \in S$ with $\|f_K\|\leq 2\|f\|$ and $f_K \res_K = f\res_K$
because the restriction map is a quotient morphism of $C^*$-algebras.
Since the net $(f_K)$ is bounded and converges to $f$ in the compact open topology,
(a) implies that it also converges in the strict topology.
Therefore $S$ is strictly dense in $C_b(X)$.
\qed

\subheadline{Tensor products of $C^*$-algebras}

Let ${\cal A}$ and ${\cal B}$ be $C^*$-algebras
and ${\cal A} \otimes {\cal B}$ their spatial $C^*$-tensor product (defined by the minimal cross norm)
([Fi96]),
which is a suitable completion of the algebraic tensor product
${\cal A} \otimes {\cal B}$, turning it into a $C^*$-algebra.
We then have homomorphisms
$$ i_{\cal A} \: M({\cal A}) \to M({\cal A} \otimes {\cal B}), \quad i_{\cal B} \: M({\cal B}) \to M({\cal A} \otimes {\cal B}), $$
uniquely determined by
$$ i_{\cal A}(\phi)(A \otimes B) = (\phi\cdot A) \otimes B, \quad
i_{\cal B}(\phi)(A \otimes B) = A \otimes (\phi\cdot B). $$
Moreover, for each complex Hilbert
space ${\cal H}$, we have
$$ \Rep({\cal A}\otimes {\cal B},{\cal H})
\cong \{(\alpha,\beta) \in \Rep({\cal A},{\cal H}) \times \Rep({\cal B},{\cal H}) \:
[\alpha({\cal A}),\beta({\cal B})] = \{0\}\}. $$
This correspondence is established by assigning to each pair
$(\alpha,\beta)$ with commuting range the representation
$$ \pi := \alpha \otimes \beta \: {\cal A} \otimes {\cal B} \to B({\cal H}), \quad
a \otimes b  \mapsto \alpha(a) \beta(b). $$
Note that this representation of ${\cal A} \otimes {\cal B}$ is non-degenerate if
$\alpha$ and $\beta$ are non-degenerate.

\Lemma A.2. The following assertions hold  for the embedding
$i_{\cal A}\: M({\cal A}) \to
M({\cal A} \otimes {\cal B})$:
\litem{(1)} The map
$$ i_{\cal A}^{-1} \: M({\cal A}) \otimes \1 \to M({\cal A}), \quad m \otimes \id_{\cal B} \mapsto m  $$
is continuous with respect to the strict topology on its domain obtained from
${\cal A} \otimes {\cal B}$  and the strict topology on its range obtained from ${\cal A}.$
\litem{(2)} Its restriction to bounded subsets is a homeomorphism.
\litem{(3)} $i_{\cal A}({\cal A})$ is dense in $M({\cal A}) \otimes \1$ with respect to
the strict topology on $M({\cal A} \otimes {\cal B})$.

\Proof.  (1) The strict topology on $M({\cal A})$ is defined by the seminorms
$$ p_a(m) = \|m \cdot a\| + \|a \cdot m\|, $$
satisfying $p_a \circ i_{\cal A}^{-1} = p_{a \otimes \1}$, which
shows immediately that $i_{\cal A}^{-1}$ is continuous.

(2) Since the embedding $i_{\cal A}$ is isometric, it suffices to show that for each
bounded subset ${\cal M} \subeq M({\cal A})$, the restriction of $i_{\cal A}$ to ${\cal M}$
is continuous. Since $i_{\cal A}$ is linear, it suffices to show that for each
bounded net $(M_\nu)$ with $\lim M_\nu = 0$ in the strict topology of
$M({\cal A})$, we also have $\lim i_{\cal A}(M_\nu) = 0$ in $M({\cal A} \otimes {\cal B})$.
For $A \in {\cal A}$ and $B \in {\cal B}$ we have
$$ \|M_\nu(A \otimes B)\| = \|M_\nu A \otimes B\| = \|M_\nu A\| \|B\| \to 0 $$
and likewise $(A \otimes B)M_\nu \to 0$. Since the elementary tensors span a dense
subset of $A \otimes B$, the boundedness of the net $(M_\nu)$ implies that
$i_{\cal A}(M_\nu)\to 0$ holds in the strict topology of $M({\cal A} \otimes {\cal B})$
(cf.\ Wegge--Olsen [WO93], Lemma~2.3.6).

(3) Let $\{E_\alpha\}$ be any approximate identity of ${\cal A}$,
satisfying $\|E_\alpha\|\leq 1$.
Then for any $A\in M({\cal A})$, the net $\{AE_\alpha\}\subset
M({\cal A})$
is bounded by $\|A\|$ and converges to $A$ in the strict
topology of $M({\cal A})$, and hence in the strict
topology of $M({\cal A}\otimes {\cal B})$ by (2). This proves~(3).
\qed

\Lemma A.3. For each non-degenerate representation $\pi \in \Rep({\cal A}
\otimes {\cal B},{\cal H})$ the representations
$\pi_1(a) := \tilde\pi(a \otimes \1)$ and
$\pi_2(b) := \tilde\pi(\1 \otimes b)$ are non-degenerate,
where $\tilde\pi$ denotes the unique extension of $\pi$ from
$\cl A.\otimes\cl B.$ to ${M(\cl A.\otimes\cl B.)}\,.$
Moreover, the corresponding
extensions
$\tilde\pi_1 \in \Rep(M({\cal A}),{\cal H})$ and
$\tilde\pi_2 \in \Rep(M({\cal B}),{\cal H})$ from $\pi_1,\;\pi_2$ on
$\cl A.,\;\cl B.$ resp., satisfy
$$\tilde\pi_1 = \tilde\pi \circ i_{{\cal A}} \quad \hbox{ and } \quad
\tilde\pi_2 = \tilde\pi \circ i_{{\cal B}}. $$
In particular, the representations
$\tilde\pi \circ i_{\cal A}$ and
$\tilde\pi \circ i_{\cal B}$ are continuous with respect to the
strict topology on $M({\cal A})$, $M({\cal B})$  resp., and the
the topology of pointwise convergence on $B({\cal H})$.

\Proof. To see that $\pi_1$ is non-degenerate, we observe that for
$a\otimes b \in {\cal A}\otimes {\cal B}$ we have
$\pi(a \otimes b) = \pi_1(a)\pi_2(b) = \pi_2(b)\pi_1(a)$, so that any
vector annihilated by $\pi_1({\cal A})$ is also annihilated by
${\cal A} \otimes {\cal B}$, hence zero. The same argument proves
non-degeneracy of~$\pi_2$.

For $m \in M({\cal A})$, we have
$$ \tilde\pi(m\otimes \1) \pi_1(a) = \tilde\pi(m\otimes \1)\tilde\pi(a \otimes \1)
= \tilde\pi(ma\otimes \1) = \pi_1(ma) = \tilde\pi_1(m)\tilde\pi_1(a), $$
so that the non-degeneracy of $\pi_1$ implies $\tilde\pi \circ i_{\cal A} = \tilde\pi_1$,
and likewise $\tilde\pi \circ i_{\cal B} = \tilde\pi_2$.

The last assertion follows from the general fact that for a non-degenerate
representation of ${\cal A}$, the corresponding extension to $M({\cal A})$
is continuous with respect to the strict topology on $M({\cal A})$ and the
topology of pointwise convergence on $B({\cal H})$; similary for ${\cal B}$.
\qed

\Lemma A.4. Let $G_1, G_2$ be topological groups and suppose that
$({\cal A}_1, \eta_1)$, resp., $({\cal A}_2, \eta_2)$ are full host algebras for
$G_1$, resp., $G_2$. Then
$$ \eta \: G_1 \times G_2 \to M({\cal A}_1 \otimes {\cal A}_2), \quad
(g_1, g_2) \mapsto i_{{\cal A}_1}(\eta_1(g_1)) i_{{\cal A}_2}(\eta_2(g_2)) $$
defines a full host algebra of $G_1 \times G_2$.

\Proof. This follows from the observation that unitary representations of
the direct product group $G := G_1\times G_2$ can be viewed as pairs
of commuting representations $\pi_j \: G_j \to U({\cal H})$,
and we have the same picture on the level of non-degenerate representations of
$C^*$-algebras.
We only have to observe that both pictures are compatible.
In fact, let $\pi_j$ be commuting unitary representations of $G_j$, $j = 1,2$,
and $\tilde\pi_j$ the corresponding representations of the host algebras ${\cal A}_j$.
Then we have
$$ \big(\eta^*(\tilde\pi_1 \otimes \tilde\pi_2)\big)(g_1, g_2)
=  (\tilde\pi_1 \otimes \tilde\pi_2)(\eta_1(g_1) \otimes \eta_2(g_2))
=  \tilde\pi_1(\eta_1(g_1)) \tilde\pi_2(\eta_2(g_2))
=  \pi_1(g_1)\pi_2(g_2).
\qeddis

Corollary~A.7 below provides a converse to this lemma.

\subheadline{Ideals of multiplier algebras}

Let ${\cal A}$ be a $C^*$-algebra and $M({\cal A})$ its multiplier algebra.
We are interested in the relation between the ideals of ${\cal A}$ and $M({\cal A})$.

\Lemma A.5. {\rm(a)} Each strictly closed ideal $J \subeq M({\cal A})$ coincides with the
strict closure of the ideal $J\cap {\cal A}$ of ${\cal A}$, which is norm-closed.

{\rm(b)} For each norm closed
ideal $I \trile {\cal A}$, its strict closure $\tilde I$ satisfies $\tilde I \cap {\cal A} = I$.

{\rm(c)} The map $J \mapsto J \cap {\cal A}$ induces a bijection from the set of strictly
closed ideals of $M({\cal A})$ onto the set of norm-closed ideals of ${\cal A}$.

\Proof. (a) Let $(u_i)_{i \in I}$ be an approximate identity in ${\cal A}$ and
$\mu \in J$. Then $\mu u_i \in J \cap {\cal A}$ converges to $\mu$ in the strict topology,
and the assertion follows.
Since on ${\cal A}$ the norm topology is finer than the strict topology, the ideal
$J \cap {\cal A}$ of ${\cal A}$ is norm-closed.

(b) The ideal $I$ is automatically $*$-invariant
([Dix64], Prop.~1.8.2), so that
${\cal A}/I$ is a $C^*$-algebra. Let $q \: {\cal A} \to {\cal A}/I$ denote the quotient
homomorphism. The existence of an approximate
identity in ${\cal A}$ implies that $I$ is invariant under the
left and right action of the multiplier algebra, so that we obtain a
natural homomorphism
$$ M(q) \: M({\cal A}) \to M({\cal A}/I), $$
which is strictly continuous ([Bu68, Prop.~3.8]).
Then $\tilde I := \ker M(q) \trile M({\cal A})$ is a strictly closed ideal satisfying
$\tilde I \cap {\cal A} = I$, and (a) implies that $\tilde I$ is the strict closure of $I$.

(c) follows from (a) and (b).
\qed

The following proposition shows that for each closed normal subgroup
$N$ of a topological group $G$ with a host algebra, the quotient group
$G/N$ also has a host algebra:

\Proposition A.6. Let $G$ be a topological group and suppose that ${\cal A}$
is a host algebra for $G$ with respect to the homomorphism
$$\eta_G \: G \to M({\cal A}). $$
Let $N \trile G$ be a closed normal subgroup,
$\tilde I_N \trile M({\cal A})$ the strictly closed ideal generated by $\eta_G(N) - \1$,
and $I_N := {\cal A} \cap \tilde I_N$.
Then $\eta_G$ factors through a homomorphism
$$ \eta_{G/N} \: G/N \to M({\cal A}/I_N), $$
turning ${\cal A}/I_N$ into a host algebra for the quotient group $G/N$.
If, in addition, ${\cal A}$ is a full host algebra of $G$,
then ${\cal A}/I_N$ is a full host algebra of $G/N$.

\Proof. If $\pi$ is a unitary representation of $G$, then we write
$\pi_{\cal A}$ for the corresponding representation of ${\cal A}$ and
$\tilde\pi_{\cal A}$ for the extension to $M({\cal A})$ with $\tilde\pi_{\cal A} \circ \eta_G = \pi$.
Further,  let $q_G \: G \to G/N$ denote the quotient map.

We consider the $C^*$-algebra ${\cal B} := {\cal A}/I_N$ and recall that
the quotient morphism
$q \: {\cal A} \to {\cal B}$ induces a strictly continuous morphism
$M(q) \: M({\cal A}) \to M({\cal B})$ ([Bu68, Prop.~3.8]).
In view of $I_N = \ker q = (\ker M(q)) \cap {\cal A}$, Lemma~A.5
implies that $\ker M(q) = \tilde I_N$.

Next we observe that
$\eta_G(N) -\id_{\cal A} \subeq \tilde I_N$ implies that $N$ acts by trivial
multipliers on the algebra ${\cal B} = {\cal A}/I_N$. We therefore obtain a
group homomorphism
$$ \eta_{G/N} \: G/N \to U(M({\cal B})) \quad \hbox{ with } \quad
\eta_{G/N} \circ q_G = M(q) \circ \eta_G.  $$

To see that $\eta_{G/N}$ turns ${\cal B}$ into a host algebra for the quotient group
$G/N$, we first note that every non-degenerate representation
$\pi \: {\cal B} \to B(H)$ can be viewed as a non-degenerate representation
$\pi_{\cal A} \: {\cal A} \to B(H)$ with $\pi_{\cal A} := \pi \circ q$. The
corresponding representations of the multiplier algebras satisfy
$$ \tilde \pi \circ M(q) = \tilde \pi_{\cal A} \: M({\cal A}) \to B(H). $$
This leads to
$$ \tilde \pi \circ \eta_{G/N} \circ q_G
= \tilde \pi \circ M(q) \circ \eta_G
= \tilde \pi_{\cal A} \circ \eta_G, $$
showing that the unitary representation  $\tilde \pi \circ \eta_{G/N}$ of
$G/N$ is continuous. We thus obtain a map
$$ \eta_{G/N}^* \: \Rep({\cal B}) \to \Rep(G/N), \quad \pi \mapsto \tilde\pi \circ
\eta_{G/N}. $$
If two representations $\pi$ and $\gamma$
of ${\cal B}$ lead to the same representation of
$G/N$, i.e.,
$$\eta_{G/N}^*(\pi)=\tilde\pi \circ \eta_{G/N}=
\tilde\gamma \circ \eta_{G/N}=\eta_{G/N}^*(\gamma), $$
then the corresponding representations of $G$ coincide, i.e.,
$\tilde \pi_{\cal A} \circ \eta_G=\tilde \gamma_{\cal A} \circ \eta_G$,
but since ${\cal A}$ is a host algebra for $G$, we have
$\pi_{\cal A}=\gamma_{\cal A}$ i.e., $\pi \circ q=\gamma\circ q$ and as
$q$ is surjective, we get $\pi=\gamma$.
%
% then the corresponding representations $\pi_i \circ q$ of ${\cal A}$ coincide,
% which implies that the corresponding representations
% $\tilde \pi_i \circ M(q) \circ \eta_G$ of $G$ coincide, which leads
% to $\pi_1 = \pi_2$.

If, in addition, $\eta_G^*$ is surjective, then
every continuous unitary representation $\pi$ of $G/N$
pulls back to a continuous unitary representation of $G$ which
defines a unique representation $\rho_{\cal A}$ of ${\cal A}$ which in turn
extends to the
representation $\tilde \rho_{\cal A}$ of $M({\cal A})$ satisfying $\tilde \rho_{\cal A} \circ
\eta_G = \pi \circ q_G$. Further,
$\tilde I_N \subeq \ker \tilde \rho_{\cal A}$ implies
$I_N \subeq \ker \rho_{\cal A}$, so that $\tilde \rho_{\cal A}$ factors
via $M(q) \: M({\cal A}) \to M({\cal B})$ through a strictly continuous
representation $\tilde \pi_{\cal B}$ of $M({\cal B})$, satisfying
$\tilde \pi_{\cal B} \circ \eta_{G/N} = \pi$. This implies that $\eta_{G/N}^*$ is also
surjective.
\qed

\Corollary A.7. Let $G_1, G_2$ be topological groups and
$G := G_1 \times G_2$. If $G$ has a full host algebra
$({\cal A},\eta)$, then $G_1$ and $G_2$ have full host algebras
$({\cal A}_1, \eta_1)$ and $({\cal A}_2, \eta_2)$ with
${\cal A} \cong {\cal A}_1 \otimes {\cal A}_2$.

\Proof. The existence of host algebras of
$G_1 \cong G/(\{\1\}\times G_2)$ and
$G_2 \cong G/(G_1 \times \{\1\})$ follows directly
from the last statement in Proposition~A.6. Now Lemma~A.4 applies.
\qed

\subheadline{Symplectic space}

\Lemma A.8. In each countably dimensional symplectic vector space
$(S,\,B)$, there exists a basis $(p_n,q_n)_{n \in \N}$ with
$$ B(p_n,q_m) = \delta_{nm} \quad \hbox{ and } \quad
B(p_n,p_m) = B(q_n,q_m) = 0 \quad \hbox{ for } \quad n,m \in \N. $$
Then $Ip_n := q_n$ and $Iq_n = -p_n$ defines a complex structure on $S$
for which
$(v,w) := B(Iv,w)$ is positive definite and hence defines a (sesquilinear)
inner product on $S$ by
$$\langle v,w\rangle:= (v,w)+iB(v,w)\,.$$
% extends to a
% positive definite hermitian form on $(S,I)$ with
% $$ B(v,w) = -\Re(Iv,w) = \Im(v,w). $$
Moreover
$\{q_n\,\mid\,n\in\N\}$ is a complex
 orthonormal basis of $S$ w.r.t. $\langle\cdot,\cdot\rangle$.
% \qeddis

\Proof. Let $(e_n)_{n \in \N}$ be a linear basis of $S$.
We construct the basis elements $p_n, q_n$ inductively as follows.
If $p_1,\ldots, p_k$ and $q_1, \ldots, q_k$ are already chosen, pick a minimal
$m$ with $e_m \not\in \span\{p_1,\ldots, p_k, q_1,\ldots, q_k\}$ and put
$$ p_{k+1} := e_m - \sum_{i=1}^k\big( B(e_m, q_i)p_i +B(p_i, e_m) q_i\big) $$
to ensure that this element is $B$-orthogonal to all previous ones.
Then pick $\ell$ minimal, such that $B(p_{k+1}, e_\ell) \not=0$, put
$$ \tilde q_{k+1} := e_\ell - \sum_{i=1}^k\big( B(e_\ell, q_i)p_i +B(p_i, e_\ell) q_i\big) $$
and pick $q_{k+1} \in \R \tilde q_{k+1}$ with
$B(p_{k+1}, q_{k+1}) = 1$.
This process can be repeated ad infinitum and produces the required
bases of $S$ because for each $k$, the span of
$p_1,\ldots, p_k, q_1,\ldots, q_k$ contains at least
$e_1,\ldots, e_k$.

That
 $\{q_n\,\mid\,n\in\N\}$ a complex orthonormal basis
w.r.t. $\langle\cdot,\cdot\rangle$
follows from the definitions.
\qed

\subheadline{Acknowledgements.}

The first author gratefully acknowledges the support of the
Sonderforschungsbereich TR12, ``Symmetries and Universality in Mesoscopic Systems''
who generously supported his visit to Germany in the Summer of 2005.
The second author wishes to express his appreciation for the
generous support he received from the Australian Research
Council for his visit to the University of New South Wales in May 2004.

\def\entries{

\[AMS93
Acerbi, F., Morchio, G., Strocchi, F.,
{\it Nonregular representations of CCR algebras and algebraic
fermion bosonization},
Proceedings of the XXV Symposium on Mathematical Physics (Tor\'un, 1992).
Rep. Math. Phys. {\bf 33} no. 1-2, 7--19 (1993)

\[Bla77  Blackadar, B., {\it Infinite tensor products of $C^*$-algebras},
Pacific J. Math. {\bf 77} (1977), 313--334

\[Bl98 Blackadar, B., ``$K$-theory for Operator Algebras,'' 2nd Ed., Cambridge Univ. Press, 1998

\[Bla06
  Blackadar, B., ``Operator Algebras,'' Encyclopaedia of Mathematical Sciences Vol. 122,
Springer-Verlag, Berlin, 2006

\[Bo74  Bourbaki, N., ``Elements of mathematics. Algebra I, Chapters 1--3''
(Reprinting of the 1974 edition) Springer-Verlag, 1992

\[BG08  Buchholz, D., Grundling, H., {\it
The resolvent algebra: A new approach to
 canonical quantum systems}, J. Funct. Anal. {\bf 254} (2008), 2725--2779

\[BR03 Bratteli, O., and D.~W.~Robinson, ``Operator Algebras and Quantum Statistical
Mechanics 1,'' 2nd ed., Texts and Monographs in Physics, Springer-Verlag, 2003

\[BR97 Bratteli, O., and D.~W.~Robinson, ``Operator Algebras and Quantum Statistical
Mechanics 2,'' 2nd ed., Texts and Monographs in Physics, Springer-Verlag, 1997
% Uniqueness and simplicity for CCR in Th. 5.2.8.

\[Bu68 Busby, R.~C., {\it Double centralizers and extensions of $C^*$-algebras},
Transactions of the Amer. Math. Soc. {\bf 132} (1968), 79--99

\[BS70 Busby, R.~C., Smith, H.~A.,  {\it Representations of Twisted Group
Algebras} Trans. Amer. Math. Soc. {\bf 149:2} (1970), 503--537

\[Dix64 Dixmier, J., ``Les $C^*$-alg\`ebres et leurs
repr\'esentations,'' Gauthier-Villars, \break Paris, 1964

\[Fi96 Fillmore, P. A., ``A User's Guide to Operator Algebras,''
Wiley, New York, 1996

\[Gl03 Gl\"ockner, H., {\it Direct limit Lie groups and manifolds},
J. Math.\ Kyoto Univ.  {\bf 43} (2003), 1--26

\[GN01 Gl\"ockner, H., and K.-H. Neeb, {\it Minimally almost periodic
abelian groups and commutative $W^*$-algebras}, in
``Nuclear Groups and Lie Groups'', Eds. E. Martin Peinador et al.,
Research and Exposition in Math. {\bf 24}, Heldermann Verlag,
2001, 163--186

\[Gr97  Grundling, H., {\it A group algebra for inductive limit groups.
Continuity problems of the canonical commutation relations},
Acta Appl. Math. {\bf 46} (1997), 107--145

\[Gr05  Grundling, H., {\it Generalizing group algebras}, J. London Math. Soc.
{\bf 72} (2005), 742--762. An erratum is in J. London Math. Soc.
{\bf 77} (2008), 270--271
% http://arxiv.org/abs/math/0410128

\[GH88
Grundling, H., Hurst, C.A., {\it  A note on regular
states and supplementary conditions}, Lett.
  Math. Phys. {\bf 15}, 205--212 (1988) [Errata: ibid. {\bf 17},
  173--174 (1989)]

\[He71 Hegerfeldt, G.C., {\it Decomposition into irreducible representations
for the Canonical Commutation Relations}, Nuovo Cimento B, {\bf 4}
(1971), 225--244

\[KR83 Kadison, R.~V., and Ringrose, J.~R., ``Fundamentals of the Theory of Operator
  Algebras II,'' New York, Academic Press 1983

\[M68 Manuceau, J., {\it C*--algebre de relations de commutation}.
Ann. Inst. Henri Poin\-car\'e {\bf 8} (1968), 139--161

\[MSTV73 Manuceau, J., Sirugue, M., Testard, D, Verbeure, A.
{\it The Smallest C*--algebra for the Canonical Commutation Relations}.
Commun. Math. Phys. {\bf 32} (1973), 231--243

\[Ne08 Neeb, K-H., {\it
A complex semigroup approach to group algebras of infinite
dimensional Lie groups}, Semigroup Forum {\bf 77} (2008), 5--35

\[PR89 Packer, J., and Raeburn, I., {\it Twisted crossed products of
$C^*$-algebras},
Math. Proc. Camb. Phil. Soc. {\bf 106} (1989) 293--311

\[Pe97 Pestov, V., {\it Abelian topological groups without irreducible
Banach representations}, in ``Abelian groups, module theory
and topology'' (Padua,
1997), 343--349, Lecture Notes in Pure and Appl. Math. {\bf 201},
Dekker, New York 1998

\[Sch90 Schaflitzel, R., {\it Decompositions of Regular Representations of the
Canonical Commutation Relations}, Publ. RIMS Kyoto Univ. {\bf 26} (1990), 1019--1047

\[Se67 Segal, I.~E., {\it Representations of the canonical commutation
relations,} Carg\'ese lectures in theoretical physics, 107--170,
Gordon and Breach, 1967

% \[Tak79 Takesaki, M., ``Theory of Operator Algebras I,''
% Berlin, Springer-Verlag, 1979

\[TZ75 Takeuti, G., Zaring, W.~M., ``Introduction to Axiomatic Set Theory,''
Berlin, Springer-Verlag, 1975

\[WO93 Wegge-Olsen, N.\ E., ``K--theory and $C^*$-algebras,''
Oxford Science Publications, 1993

\[Wo95 Woronowicz, S. L., {\it $C^*$-algebras generated by unbounded elements},
Rev. Math. Phys. {\bf 7} (1995), 481--521

}

\references

\vfill\eject
\bye

\maketitle
% \begin{abstract}
% \end{abstract}
%%%%%%%%%%%%%%%%%%%%%%%%%%%%%%%%%%%%%%%%%%%%%%%%%%%%%%%%%%%%%%%%%%%%%%%%%%%%%%
%\section*{}

The proof of the main theorem of the paper~\cite{GrNe} contains an error.
We are grateful to Prof.~Ralf Meyer (Mathematisches Institut, Georg-August Universit\"at G\"ottingen) for pointing
out this mistake.

\section{Review}

The example constructed in~\cite{GrNe}, was intended to be an example
of a full host algebra for a  non--locally compact group, in our project of generalising group algebras to
general topological groups.
The concept of a group algebra is generalized as follows:
\begin{defi}  Let $G$ be a topological group.
A {\it host algebra for $G$} is a pair
$({\cal L}, \eta)$, where  ${\cal L}$ is a $C^*$-algebra and
$\eta \colon G \to \U(M({\cal L}))$ is a group homomorphism
such that:
\begin{itemize}
\item[\rm(H1)] For each non-degenerate representation $(\pi, {\cal H})$
of $\cL,$ the representation $\widetilde \pi \circ \eta$ of $G$ is continuous,
where $\tilde\pi$ denotes the canonical extension of $\pi$ to the
multiplier algebra $M(\cL).$
\item[\rm(H2)] For each complex Hilbert space
${\cal H},$ the corresponding map
$$ \eta^* \colon \rep({\cal L},{\cal H}) \to \rep(G, {\cal H}), \quad
\pi \mapsto \widetilde \pi \circ \eta $$
is injective. We write $\rep(G,{\cal H})_\eta$ for the range of $\eta^*.$
\end{itemize}
%and its elements are called {\it $\cL$-representations of $G$} on $\cH$.
%Note that $\eta^*$ depends on $\cal H$.
% \\[1mm]
We call $({\cal L}, \eta)$  a {\it full host algebra} if, in addition, we have:
\begin{itemize}
\item[\rm(H3)] $\rep(G,{\cal H})_\eta = \rep(G,{\cal H})$
for each Hilbert space~${\cal H}.$
\end{itemize}
\end{defi}
Thus by (H2) and (H3), a full host algebra, when it exists, carries precisely the continuous
unitary representation theory of $G,$ and if it is not full, it carries some subtheory of the
continuous unitary representations of $G$. In \cite[Prop.~3.4]{GrNe} we see that there is a large
collection of (non--full) host algebras for the context of  \cite{GrNe}, in fact enough to include
all cyclic representations by Proposition~\ref{cyclicrep} below. Further examples of host algebras
for a non--locally compact group, can be found
in~\cite{NeSaZe}, where the representation subtheories are given by  a spectral condition.
 However, a {\it full} host algebra
is hard to construct in the  non--locally compact context, and given the failure of the example
in \cite{GrNe}, the construction of one remains  an open problem.

\section{Counterexample}

Recall the context of \cite{GrNe}.
Let $(S,B)$ be a countably infinite dimensional symplectic vector space
and $S_\sigma = S \times \T$ be the corresponding Heisenberg group with
the product $(v,z)(v',z') = (v + v', zz' \sigma(v,v'))$ where
$\sigma$ is the 2-cocycle $\sigma(\cdot,\cdot)=\exp(i B(\cdot,\cdot)/2).$
Fix a symplectic basis $(q_j,p_j)_{j\in\N}$ with finite span space $S.$
For each symplectic subspace $S_j:={{\rm Span}(q_j,p_j)}$ construct
a twisted convolution C*-algebra ${\cal L}_j$
w.r.t. $\sigma.$
This is isomorphic to
the compact operators $ \cK(\ell^2(\N)).$
Fix a countable approximate identity of projections $(E_n^{(j)})_{n \in \N}$ in
${\cal L}_j$ for each $j.$
We associate to each sequence
$\boldn = (n_1, n_2, \ldots)$ of natural numbers
an infinite tensor product algebra $\cL[\boldn]$ of all ${\cal L}_j$ in the sense of Blackadar, using the
sequence of projections $(E_{n_j}^{(j)})_{j \in \N}.$
Every such algebra acts
faithfully by multipliers on $\cL[\bold1]$ where $\bold1=(1,1,\ldots)$ and we write
$\cL[E] \subseteq M(\cL[\bold1])$ for the subalgebra generated by all these
infinite  tensor products
$\cL[\boldn]$. The Weyl algebra $\cal A$ w.r.t. $(S,B)$ acts as  multipliers on
each  $\cL[\boldn],$ compatibly w.r.t. the products between the different $\cL[\boldn],$ hence
we have an embedding ${\cal A}\subset M({\cal L}[E]).$ Denote the generating unitaries of
 $\cal A$ by $\delta_s,$ $s\in S.$ The regular representations of  $\cal A$ on $\cl H.$
 correspond to ${\rm Rep}((S,\sigma),\cl H.).$

In \cite[Theorem~3.6]{GrNe}, we defined the monomorphism
$\eta:S_\sigma\to U(M({\cal L}[E]))$  by
$$\eta((s,t)):= t\delta_s\in {\cal A}\subset M({\cal L}[E]),$$
and claimed that
${\cal L}[E]$ is a host algebra, i.e. the map
$$\eta^*:{\rm Rep}({\cal L}[E],\cl H.)\to {\rm Rep}((S,\sigma),\cl H.)$$
is injective.
The proof of this, given in the second paragraph of the proof, is invalid,
as the following counterexample to the claim demonstrates.

\begin{eje}
(R. Meyer)\\
We will construct two distinct representations in ${\rm Rep}({\cal L}[E],\cl H.)$
which map to the same image with respect to  $\eta^*.$ Thus $\eta^*$ is not injective.

Since $S$ is countably infinite dimensional, ${\cal L}[E]$
has the following proper nontrivial ideals: ${\cal L}[{\bf 1}]$ and %(cf. \cite[Prop.~3.5(ii)]{GrNe}),
\[
\cl I.[\b n.]=\big[\big\{ {\cal L}[{\bf q}]_0 \,\mid\, [\b q.]\leq[\b n.]\big\}\big]
\]
(cf. \cite[Prop.~3.5(iii)]{GrNe}) for any sequence
${\bf n}\in\N^\infty$ 
where the square brackets enclosing the set on the right, indicate closed span.
Here we let $[\b n.]$ denote
 the equivalence class of ${\bf n}$ in
$\N^\infty\big/\mathord{\sim}$ where ${\bf n}\sim{\bf m}$ if they
differ only in finitely many entries.
We denote $[\b n.]\geq [\b m.]$ if
for any representatives $\b n.$ and $\b m.$ resp., we have that
there is an $N$ (depending on the representatives)
such that $n_k\geq m_k$ for all $k>N\,.$

Fix $[{\bf n}]\not=[{\bf 1}]$,
then we can also define the closed ideal
\[
\cl I.[\b n.]^o:=\big[\big\{ {\cal L}[{\bf q}]_0 \,\mid\, [\b q.]<[\b n.],\;\hbox{and}\; [\b q.]\not=[\b n.] \big\}\big]
\subset\cl I.[\b n.] \,.
\]
The inclusion is proper, because by the proof of \cite[Prop.~3.5(ii)]{GrNe}
we have
% \begin{equation}
%  \label{eq:e1}
\[
{\cal I}[{\bf n}]^o\cap{\cal L}[{\bf n}]=\{0\}.
\]
%\end{equation}
Denote the quotient map by $\zeta:\cl I.[\b n.]\to\cl I.[\b n.]\big/\cl I.[\b n.]^o.$
Thus $\zeta$ is faithful on ${\cal L}[{\bf n}]$, and as $\zeta$ is zero on each ${\cal L}[{\bf q}]$
for $ [\b q.]<[\b n.]$ and $[\b q.]\not=[\b n.],$ and respects
$C^*$-operations, we have
\[
\zeta\big(\cl I.[\b n.]\big)=\zeta\big(\cl L.[\b n.]\big)\cong\cl L.[\b n.].
\]
We conclude that
%$\cl I.[\b n.]^o+\cl L.[\b n.]\subseteq\cl I.[\b n.] $ is closed by \cite[Thm.~2.1]{ZhDi},
%hence by denseness, that
$\cl I.[\b n.]^o+\cl L.[\b n.]=\cl I.[\b n.].$
By ${\cal I}[{\bf n}]^o\cap{\cal L}[{\bf n}]=\{0\}$
we have  $\cl I.[\b n.]=\cl I.[\b n.]^o\oplus\cl L.[\b n.].$

Let $\pi_1\in{\rm Rep}({\cal L}[E],\cl H.)$ be a representation which is nondegenerate on
${\cal L}[{\bf 1}].$
%(hence it is the unique restriction of the corresponding extension to $M(\cL[\bold1])$).
Then $\pi_1$ is nondegenerate on each ${\cal L}[{\bf n}]$ because
\[
\pi_1({\cal L}[{\bf n}])\cl H.=\pi_1({\cal L}[{\bf n}])\big(\pi_1({\cal L}[{\bf 1}])\cl H.\big)
=\pi_1({\cal L}[{\bf 1}])\cl H.=\cl H..
\]
Define a representation $\widehat\pi_1\in{\rm Rep}(\cl I.[\b n.],\cl H.)$ by
\[
\widehat\pi_1(K+L):=\pi_1(L)\quad\hbox{for}\quad K\in\cl I.[\b n.]^o\;\hbox{and}\;L\in\cl L.[\b n.]
\]
using the decomposition $\cl I.[\b n.]=\cl I.[\b n.]^o\oplus\cl L.[\b n.].$
Let  $\pi_2\in{\rm Rep}({\cal L}[E],\cl H.)$ be the unique (strict) extension of $\widehat\pi_1$
from $\cl I.[\b n.]$ to ${\cal L}[E].$ Clearly $\pi_1\not=\pi_2$ as $\pi_1({\cal L}[{\bf 1}])
\not=\{0\}=\pi_2({\cal L}[{\bf 1}]),$ but $\pi_1\restriction{\cal L}[{\bf n}]
=\pi_2\restriction{\cal L}[{\bf n}]$ (hence $\pi_2$ is nondegenerate).

Recall that for a representation $\pi\in{\rm Rep}({\cal L}[E],\cl H.),$ its image
$\eta^*\pi\in {\rm Rep}((S,\sigma),\cl H.)$ is uniquely determined by
\[
(\eta^*\pi)(s)\cdot \pi(L)\psi=\pi\big(\eta(s)L)\psi
\quad \mbox{ for all } \quad s\in S, L\in{\cal L}[E]\quad \mbox{ and } \quad
\psi\in\cl H..\]  By this formula, it is enough to
let $L$ range over a subset $R\subset{\cal L}[E]$ such that  $\pi(R)\cl H.$ is dense in $\cl H..$
For $\pi_1$ and $\pi_2$ above, such a set is ${\cal L}[{\bf n}],$ as $\pi_1$ and $\pi_2$ restricted to
${\cal L}[{\bf n}]$ are nondegenerate. Moreover, as these restrictions to ${\cal L}[{\bf n}]$ coincide,
it follows that $\eta^*\pi_1=\eta^*\pi_2.$ Thus  $\eta^*$ is not injective.
\end{eje}

\section{Surviving results}

We discuss what remains of the paper after this mistake.
We do not know whether $(S,\sigma)$ has a full host algebra.
Whilst we lose the full host algebra of Theorems~3.1 and 3.6, we still have the family of
host algebras ${\cal L}[{\bf n}]$, indexed by the equivalence classes
$[{\bf n}]\in\N^\infty\big/\mathord{\sim}$. By the proposition below, for each (regular) cyclic representation
$\pi\in {\rm Rep}((S,\sigma),\cl H.)$ there is an $[{\bf n}]\in\N^\infty\big/\mathord{\sim}$
and a $\pi_0\in {\rm Rep}({\cal L}[{\bf n}],\cl H.)$ such that $\pi=\eta^*\pi_0.$
Thus we obtain a structure theory for ${\rm Rep}((S,\sigma),\cl H.).$
To be precise, denote
\[
 {\rm Rep}_{[{\bf n}]}((S,\sigma),\cl H.):=\eta^*\Big({\rm Rep}({\cal L}[{\bf n}],\cl H.)\Big),
\]
and recall from \cite[Prop.~3.4]{GrNe} that
\[
 {\rm Rep}_{[{\bf n}]}((S,\sigma),\cl H.)=
 \Big\{\pi\in {{\rm Rep}\big((S,\sigma),\cl H.\big)}\;\mid\;\slim\limits_{k\to\infty}\,\slim_{m\to\infty}
 \pi(E_{n_k}^{(k)})\cdots\pi(E_{n_m}^{(m)})=\un\Big\}.
\]
From this it follows that if $[{\bf n}]\leq[{\bf p}]$ then ${\rm Rep}_{[{\bf n}]}((S,\sigma),\cl H.)\subseteq
{\rm Rep}_{[{\bf p}]}((S,\sigma),\cl H.)$  as $E_{n_k}^{(k)}\leq E_{p_k}^{(k)}$ for all $k$ except a finite number.
Thus, for a set of representations
\[
\{\pi_j\in {\rm Rep}_{[{\bf n}_j]}((S,\sigma),\cl H._j)\,\mid\, j\in J\}
\]
for an index set $J,$ which is bounded above in the sense that there is a $[{\bf k}]\in\N^\infty\big/\mathord{\sim}$ such that
$[{\bf n}_j]\leq[{\bf k}]$ for all $j\in J,$ then by  $\pi_j\in {\rm Rep}_{[{\bf k}]}((S,\sigma),\cl H._j),$
it follows for the direct sum that
\[
\bigoplus_{j\in J}\pi_j\in {\rm Rep}_{[{\bf k}]}
\Big((S,\sigma),\mathop{\oplus}\limits_{j\in J}\cl H._j\Big).
\]
Thus each finite direct sum of cyclic representations is in some ${\rm Rep}_{[{\bf k }]}((S,\sigma),\cl H.).$

Conversely, each $\pi\in{\rm Rep}_{[{\bf n}]}((S,\sigma),\cl H.)$ for separable $\cl H.$ has a
direct integral decomposition in terms of irreducibles in
${\rm Rep}_{[{\bf n}]}((S,\sigma),\cl H.')$ for some $\cl H.'$ (consider the argument below the proof of
\cite[Theorem~3.6]{GrNe}).
\begin{pro}
\label{cyclicrep}
The directed union
$\displaystyle{\bigcup_{[{\bf n}]\in\N^\infty\big/\mathord{\sim}}
{\rm Rep}_{[{\bf n}]}((S,\sigma),\cl H.)}$ contains all cyclic representations in
${\rm Rep}((S,\sigma),\cl H.).$
\end{pro}
\begin{beweis}
Consider a cyclic representation $\pi\in{\rm Rep}((S,\sigma),\cl H.)$ with a cyclic unit vector $\varphi\in \cl H..$
We need to prove that there is an ${\bf n}\in\N^\infty$ such that
\[
\slim\limits_{k\to\infty}\,\slim_{m\to\infty}\pi(E_{n_k}^{(k)})\cdots\pi(E_{n_m}^{(m)})=\un
 \]
As $\pi$ is nondegenerate on each $(S_j,\sigma),$ its extension to $\cl L._j$ is nondegenerate (still denoted by~$\pi$).
As $(E_{\ell}^{(j)})_{\ell\in\N}\subset \cl L._j$ is an approximate identity for $\cl L._j,$ it follows that
$\slim\limits_{\ell\to\infty}\pi(E_{\ell}^{(j)})=\un$ for all $j\in\N.$ Given an ${\bf n}\in\N^\infty,$ abbreviate
$P_k:=\pi(E_{n_k}^{(k)})$ for all $k,$ then using
\[
P_1\cdots P_m-\un=P_2\cdots P_m(P_1-\un)+P_3\cdots P_m(P_2-\un)+\cdots+P_m(P_{m-1}-1)+(P_m-\un)
\]
we conclude that
\[
\big\|(P_1\cdots P_m-\un)\varphi\big\|\leq \sum_{i=1}^m\|(P_i-\un)\varphi\|.
\]
By $\slim\limits_{\ell\to\infty}\pi(E_{\ell}^{(j)})=\un$ it is possible to find an ${\bf n}\in\N^\infty$ for
which the series
$
\sum\limits_{i=1}^\infty\|(P_i-\un)\varphi\|
$
converges. Fix such a choice of ${\bf n}.$ Thus by
\[
\big\|\big(\pi(E_{n_k}^{(k)})\cdots\pi(E_{n_m}^{(m)})-\un\big)\varphi\big\|\leq
\sum\limits_{i=k}^\infty\|(P_i-\un)\varphi\|
\]
we conclude that
\[
\lim_{k\to\infty}\,\lim_{m\to\infty}\big\|\big(\pi(E_{n_k}^{(k)})\cdots\pi(E_{n_m}^{(m)})
-\un\big)\psi\big\|=0 \quad \mbox{ for } \quad
\psi =\varphi.\]
We need to prove this relation  for any $\psi$ in a dense subspace of $\cl H.,$ where we keep ${\bf n}$  fixed.
As $\varphi$ is cyclic, consider the dense space ${\rm Span}\,\pi(S)\varphi$. As it consists of finite linear combinations we only
need to show that
\[
\lim_{k\to\infty}\,\lim_{m\to\infty}\big\|\big(\pi(E_{n_k}^{(k)})\cdots\pi(E_{n_m}^{(m)})-\un\big)\pi(s)\varphi\big\|=0
\]
for all $s\in S$ for the fixed ${\bf n}$ above. As $S$ consists of finite real combinations of the given
symplectic basis $(q_j,p_j)_{j\in\N}$, for each $s\in S,$ $\pi(s)$ will commute with all $\pi(E_{\ell}^{(k)})$ for $k$ large enough.
Hence
\begin{eqnarray*}
\lim_{k\to\infty}\,\lim_{m\to\infty}\big\|\big(\pi(E_{n_k}^{(k)})\cdots\pi(E_{n_m}^{(m)})\!\!&-&\!\!\un\big)\pi(s)\varphi\big\|
\\[1mm]
&=&\lim_{k\to\infty}\,\lim_{m\to\infty}\big\|\pi(s)\big(\pi(E_{n_k}^{(k)})\cdots\pi(E_{n_m}^{(m)})-\un\big)\varphi\big\|\\[1mm]
&=&\lim_{k\to\infty}\,\lim_{m\to\infty}\big\|\big(\pi(E_{n_k}^{(k)})\cdots\pi(E_{n_m}^{(m)})-\un\big)\varphi\big\|=0\,.
\end{eqnarray*}
This proves the proposition.
\end{beweis}

% \section*{Acknowledgements.}

\providecommand{\bysame}{\leavevmode\hbox to3em{\hrulefill}\thinspace}